\documentclass[aihp]{imsart}

\RequirePackage{amsthm,amsmath,amsfonts,amssymb}
\RequirePackage[numbers]{natbib}
\RequirePackage[colorlinks,citecolor=blue,urlcolor=blue]{hyperref}
\RequirePackage{graphicx}

\startlocaldefs
\theoremstyle{plain}

 \newtheorem{thm0}{Theorem}[section]
 \newtheorem{con0}{Theorem}[section]
 \newtheorem{exa0}{Theorem}[section]

\theoremstyle{remark}

 \newtheorem{con1}[con0]{\it{Condition}}
 \newtheorem{lem1}[thm0]{Lemma}
 \newtheorem{thm1}[thm0]{Theorem}
 \newtheorem{cor1}[thm0]{Corollary}
 \newtheorem{pro1}[thm0]{Proposition}
 \newtheorem{rem1}[thm0]{Remark}
 \newtheorem{exa1}[exa0]{\it{Example}}

 \usepackage{color}
 \definecolor{MyDarkBlue}{rgb}{0,0.08,0.45}
 \definecolor{darkgreen}{rgb}{0.00,0.50,0.00}

 \def\bcondition{\begin{con1}}\def\econdition{\end{con1}}
 \def\blemma{\begin{lem1}}\def\elemma{\end{lem1}}
 \def\btheorem{\begin{thm1}}\def\etheorem{\end{thm1}}
 \def\bcorollary{\begin{cor1}}\def\ecorollary{\end{cor1}}
 \def\bproposition{\begin{pro1}}\def\eproposition{\end{pro1}}
 \def\bremark{\begin{rem1}}\def\eremark{\end{rem1}}
 \def\bexample{\begin{exa1}}\def\eexample{\end{exa1}}

 \def\benumerate{\begin{enumerate}}\def\eenumerate{\end{enumerate}}
 \def\bitemize{\begin{itemize}}\def\eitemize{\end{itemize}}\def\itm{\item}

 \def\beqlb{\begin{eqnarray}}\def\eeqlb{\end{eqnarray}}
 \def\beqnn{\begin{eqnarray*}}\def\eeqnn{\end{eqnarray*}}

 \def\eqref#1{{\rm(\ref{#1})}}

 \def\ar{&}\def\nnm{\nonumber}
 \def\qqquad{\qquad\qquad}

 \def\qed{\hfill$\square$\smallskip}

 \def\mrm{\mathrm}\def\mbf{\mathbf}\def\mcr{\mathcal}
 \def\mbb{\mathbb}

 \def\d{\mrm{d}}\def\e{\mrm{e}}
 
 \def\I{\mbf{1}}

 \def\itDelta{{\it\Delta}}\def\itPhi{{\it\Phi}}
 \def\itGamma{{\it\Gamma}}\def\itPsi{{\it\Psi}}
 \def\itOmega{{\it\Omega}}

 \def\supp{\mrm{supp}}

\def\rv{}\def\rrv{}

\endlocaldefs

\begin{document}

\begin{frontmatter}
\title{Exponential ergodicity of branching processes \\ with immigration and competition}
\runtitle{Exponential ergodicity of branching processes}

\begin{aug}
\author[A]{\inits{P.-S.}\fnms{Pei-Sen}~\snm{Li}\ead[label=e1]{peisenli@bit.edu.cn}\orcid{0000-0002-3949-9529}},
\author[B]{\inits{Z.}\fnms{Zenghu}~\snm{Li}\ead[label=e2]{lizh@bnu.edu.cn}\orcid{0000-0003-3641-4400}}
\author[C]{\inits{J.}\fnms{Jian}~\snm{Wang}\ead[label=e3]{jianwang@fjnu.edu.cn}\orcid{0000-0002-3870-4463}}
\and
\author[D]{\inits{Z.}\fnms{Xiaowen}~\snm{Zhou}\ead[label=e4]{xiaowen.zhou@concordia.ca}\orcid{0000-0002-8205-0217}}

\address[A]{School of Mathematics and Statistics, Beijing Institute of Technology, Beijing, China\printead[presep={,\ }]{e1}}

\address[B]{School of Mathematical Sciences, Beijing Normal University, Beijing, China\printead[presep={,\ }]{e2}}

\address[C]{School of Mathematics and Statistics, Fujian Normal University, Fuzhou, China\printead[presep={,\ }]{e3}}

\address[D]{Department of Mathematics and Statistics, Concordia University,  Montreal, Canada\printead[presep={,\ }]{e4}}
\end{aug}

\begin{abstract}
We study the ergodic property of a continuous-state branching process with immigration and competition. The exponential ergodicity in a weighted total variation distance is proved under natural assumptions. The main theorem applies to subcritical, critical and supercritical branching mechanisms, {\rrv including all those of} stable types. The proof is based on the construction of a Markov coupling process and the choice of a nonsymmetric control function for the distance. Those are designed to identify and to take the advantage of the dominating factor from the branching, immigration and competition mechanisms in different parts of the state space. {\rv The approach provides a way of finding a lower bound of the ergodicity rate.}
\end{abstract}

\begin{abstract}[language=french]
Nous \'etudions la propri\'et\'e ergodique d'un processus de branchement en temps et espace continu avec l'immigration et la comp\'etition. L'ergodicit\'e exponentielle dans une distance de variation totale pond\'er\'ee est prouv\'ee sous des hypoth\`eses naturelles. Le th\'eor\`eme principal s'applique aux m\'ecanismes de branchement sous-critiques, critiques et sur-critiques, y compris tous les types stables. La d\'emonstration est bas\'ee sur la construction d'un processus Markovien de couplage et le choix d'une fonction de contr\^ole non sym\'etrique pour la distance. Ceux-ci sont con\c{c}us pour identifier et profiter du facteur dominant des m\'ecanismes de branchement, d'immigration et de comp\'etition dans les parties diff\'erentes de l'espace d'\'etats.
{\rv Cette approche permet de trouver une borne inf\'erieure de la vitesse d'ergodicit\'e. }
\end{abstract}



\end{frontmatter}

\section{Introduction}

Classical \textit{Galton--Watson branching processes} are Markov processes taking values of nonnegative integers. They are models for the evolution of populations where the progenies of individuals are described by i.i.d.\ random variables. Standard references on those processes are Athreya and Ney \cite{AtN72} and Harris \cite{Har63}. {\rv The study} of \textit{continuous-state branching processes} (CB-processes) was initiated by Feller \cite{Fel51}, who noticed that a one-dimensional diffusion process may arise as the limit of a sequence of rescaled Galton--Watson branching processes. The result was extended by Lamperti \cite{Lam67a} to the situation where the limiting process may have discontinuous sample paths; see also Aliev and Shchurenkov \cite{AlS82} and Grimvall \cite{Gri74}. Let $\itPsi$ be a function on $[0,\infty)$ with the L\'evy--Khintchine representation:
 \beqlb\label{eq1.1}
\itPsi(\lambda)= b\lambda + c\lambda^2 + \int_0^\infty \big(\e^{-\lambda z}-1+\lambda z\I_{\{z\le 1\}}\big)\mu(\d z), \quad \lambda\ge 0,
 \eeqlb
where $b\in \mbb{R}$ and $c\ge 0$ are constants and $(1\land z^2)\mu(\d z)$ is a finite measure on $(0,\infty)$. The transition semigroup $(Q_t)_{t\ge 0}$ of a CB-process with \textit{branching mechanism} $\itPsi$ is defined by
 \beqlb\label{eq1.2}
\int_{\mbb{R}_+} \e^{-\lambda y} Q_t(x,\d y)= \e^{-xv_t(\lambda)}, \quad x\ge 0,\lambda> 0,
 \eeqlb
where $t\mapsto v_t(\lambda)$ is the unique strictly positive solution to the differential equation
 \beqlb\label{eq1.3}
\frac{\partial}{\partial t} v_t(\lambda)= -\itPsi(v_t(\lambda)), \quad v_0(\lambda)= \lambda.
 \eeqlb
The CB-process is eventually degenerate in the sense that it tends to either zero or infinity as $t\to \infty$. In fact, the process may explode at a finite time {\rrv with strictly positive probability. It almost surely has an infinite lifetime if and only if its branching mechanism satisfies the following \textit{conservativeness condition}:}
 \beqlb\label{eq1.4}
\int_{0+}\frac{\d\lambda}{0\vee[-\itPsi(\lambda)]}= \infty;
 \eeqlb
see, {\rrv e.g., Grey \cite[p.670]{Gre74}.} In order that the integral on the left-hand side of \eqref{eq1.4} is convergent, we have necessarily that $\itPsi'(0)= -\infty$. Let $v_t(0)= \lim_{\lambda\to 0+}v_t(\lambda)$ for $t\ge 0$.
Under condition \eqref{eq1.4}, we have $v_t(0)= 0$ for every $t\ge 0$, which is also the unique solution to \eqref{eq1.3} with $\lambda= 0$. The first moments of the transition probabilities of the CB-process are given by
 \beqlb\label{eq1.5}
\int_{\mbb{R}_+} y Q_t(x,\d y)= x\exp\big\{-\itPsi'(0+)t\big\}, \quad t,x> 0,
 \eeqlb
where
 \beqnn
\itPsi'(0+)= b-\int_1^\infty z \mu(\d z).
 \eeqnn
The branching mechanism $\itPsi$ is said to be \textit{subcritical}, \textit{critical} or \textit{supercritical} according as $\itPsi'(0+)>0$, $\itPsi'(0+)=0$ or $\itPsi'(0+)<0$, respectively. Let $\sigma\ge 0$, $a\in \mbb{R}$ and $0<\alpha< 2$ be constants. The \textit{stable branching mechanism} is defined by
 \beqlb\label{eq1.6}
\itPsi(\lambda)= \left\{\begin{array}{ll}
 a\lambda + c\lambda^2 + \sigma\lambda^\alpha, &\quad 1<\alpha< 2, \cr
 a\lambda + c\lambda^2 + \sigma\lambda\log\lambda, &\quad \alpha= 1, \cr
 a\lambda + c\lambda^2 - \sigma\lambda^\alpha, &\quad 0<\alpha< 1.
\end{array}\right.
 \eeqlb

Clearly, the conservativeness condition \eqref{eq1.4} fails for the stable branching mechanism if $\sigma>0$ and $0<\alpha< 1$. If $\sigma>0$ and $0<\alpha\le 1$, then $\itPsi'(0+)= -\infty$ and the CB-process has {\rrv infinite first moments} by \eqref{eq1.5}.

{\rv

\bexample\label{ex1.0a1} Suppose that $\itPsi(\lambda)= a\lambda + \sigma\lambda^{1+\alpha}$, where $\sigma\ge 0$, $a\in \mbb{R}$ and $1<\alpha\le 2$. Then $\itPsi'(0)= a> -\infty$ and condition \eqref{eq1.4} is satisfied. By solving \eqref{eq1.3} in this case, we see that
 \beqnn
v_t(\lambda)
 =
\frac{\e^{-at}\lambda}{\big[1 + \sigma q_\alpha(a,t)\lambda^{\alpha-1}\big]^{1/(\alpha-1)}}, \quad t\ge 0, \lambda\ge 0,
 \eeqnn
where $q_\alpha(a,t)= a^{-1}(1-\e^{-(\alpha-1)at})$ for $a\neq 0$ and $q_\alpha(0,t)= (\alpha-1)t$. \eexample

\bexample\label{ex1.0a2} Suppose that $\itPsi(\lambda)= a\lambda + \sigma\lambda\log\lambda$, where $\sigma\ge 0$ and $a\in \mbb{R}$. Then condition \eqref{eq1.4} is satisfied. In this case, we have
 \beqnn
v_t(\lambda)
 =
\exp\big\{\e^{-\sigma t}\log\lambda - a\rho(\sigma,t)\big\}, \quad t\ge 0, \lambda\ge 0,
 \eeqnn
where $\rho(\sigma,t)= \sigma^{-1}(1-\e^{-\sigma t})$ for $\sigma> 0$ and $\rho(0,t)= t$. In particular, if $\sigma> 0$, then
 \beqnn
\lim_{t\to \infty} v_t(\lambda)
 =
v_\infty:= \e^{-a/\sigma}, \quad \lambda> 0,
 \eeqnn
which implies that, by weak convergence on $[0,\infty]$,
 \beqnn
\lim_{t\to \infty} Q_t(x,\cdot)
 =
\e^{-xv_\infty}\delta_0 + (1-\e^{-xv_\infty})\delta_\infty, \quad x\ge 0.
 \eeqnn
\eexample

\bexample\label{ex1.0a3} Suppose that $\itPsi(\lambda)= a\lambda + \sigma\lambda^\alpha$, where $\sigma\ge 0$, $a\in \mbb{R}$ and $0<\alpha< 1$. Then condition \eqref{eq1.4} is not satisfied. By solving \eqref{eq1.3} in this case, we have
 \beqnn
v_t(\lambda)
 =
\big[\sigma p_\alpha(a,t)+\lambda^{1-\alpha}\e^{-(1-\alpha)at}\big]^{1/(1-\alpha)}, \quad t\ge 0, \lambda> 0,
 \eeqnn
where $p_\alpha(a,t)= a^{-1}(1-\e^{-(1-\alpha)at})$ for $a\neq 0$ and $p_\alpha(0,t) = (1-\alpha)t$. In particular, if $a> 0$, then
 \beqnn
v_t(0):= \lim_{\lambda\to 0} v_t(\lambda)
 =
[\sigma p_\alpha(a,t)]^{1/(1-\alpha)}, \quad t\ge 0
 \eeqnn
and
 \beqnn
\lim_{t\to \infty} v_t(\lambda)
 =
v_\infty:= \Big(\frac{\sigma}{a}\Big)^{1/(1-\alpha)}, \quad \lambda\ge 0.
 \eeqnn
The last limit implies that, by weak convergence on $[0,\infty]$,
 \beqnn
\lim_{t\to \infty} Q_t(x,\cdot)
 =
\e^{-xv_\infty}\delta_0 + (1-\e^{-xv_\infty})\delta_\infty.
 \eeqnn
\eexample

}

The CB-processes defined by \eqref{eq1.2} and \eqref{eq1.3} involve rich mathematical structures and have been applied to the research in several important areas. In particular, Bertoin and Le~Gall \cite{BeL00} gave a representation of the genealogical structure of the general CB-processes by Bochner's subordination and used the representation to give a precise description of the connection between the so-called \textit{Neveu's CB-process} with branching mechanism $\lambda\mapsto \lambda\log\lambda$ and the coalescent processes introduced by Bolthausen and Sznitman \cite{BoS98} in the study of spin glasses. Berestycki et al.\ \cite{BBS13} proved Neveu's CB-process may arise in a limit theorem of some rescaled particle systems. General stochastic flows associated with coalescent processes and CB-processes were studied by Bertoin and Le~Gall \cite{BeL03, BeL05, BeL06}. Constructions of the flows by strong solutions of stochastic equations were given by Dawson and Li \cite{DaL12}.

\textit{Continuous-state branching processes with immigration} (CBI-processes) were introduced by Kawazu and Watanabe \cite{KaW71} in their study of scaling limits of discrete branching processes with immigration; see also Aliev \cite{Ali85}. Let $\itPhi$ be an \textit{immigration mechanism}, which is a function on $[0,\infty)$ with the L\'evy--Khintchine representation:
 \beqlb\label{eq1.7}
\itPhi(\lambda)= \beta\lambda + \int_0^\infty \big(1-\e^{-z\lambda}\big)\nu(\d z), \quad \lambda\ge 0,
 \eeqlb
where $\beta\ge 0$ and $(1\land z)\nu(\d z)$ is a finite measure on $(0,\infty)$. {\rv Then a CBI-process with parameters $(\itPsi,\itPhi)$ {\rrv has transition semigroup $(Q^{\itPhi}_t)_{t\ge 0}$} characterized} by
 \beqlb\label{eq1.8}
\int_{\mbb{R}_+} \e^{-\lambda y} Q^{\itPhi}_t(x,\d y)= \exp\bigg\{-xv_t(\lambda) - \int_0^t \itPhi(v_s(\lambda))\d s\bigg\}, \quad \lambda> 0,
 \eeqlb
where $t\mapsto v_t(\lambda)$ is defined by \eqref{eq1.3}. The CBI-process is a natural model for the study of ergodicity, where the immigration ensures that the stationary distribution is not the degenerate distribution at zero. A sufficient and necessary integrability condition for the ergodicity in weak convergence of the process was announced in Pinsky \cite{Pin72}; see Li \cite{Li22} for a proof of the result. The exponential ergodicity of the CBI-process in the total variation distance was proved in Li and Ma \cite{LiM15} and the corresponding results for measure-valued {\rv processes were established in Friesen et al.\ \cite{Fri23} and Li \cite{Li21},} which improve the earlier results of Stannat \cite{Sta03a, Sta03b} for processes with special branching mechanisms. The exponential ergodicity in a suitably chosen Wasserstein distance was established by Friesen et al.\ \cite{FJR20} for affine Markov processes, which include finite-dimensional CBI-processes as a special case. The key tools of those explorations are the characterizations of the transition probabilities by Laplace transforms given as in \eqref{eq1.3} and \eqref{eq1.8}. Since typical critical or supercritical CBI-processes eventually go to infinity, the studies of their ergodicities have mainly been focused on subcritical mechanisms. To get the ergodicity in the total variation distance, Li and Ma \cite{LiM15} also assumed the following \textit{Grey's condition}: {\rv $\itPsi(\lambda)> 0$ for sufficiently large $\lambda$ and}
 \beqlb\label{eq1.9}
\int^{\infty-}\frac{\d\lambda}{\itPsi(\lambda)}< \infty.
 \eeqlb
Those restrictions unfortunately exclude {\rv some interesting branching mechanisms as those in Examples~\ref{ex1.0a2} and~\ref{ex1.0a3}.}

In this work we provide a general framework for exponential ergodicity applicable to branching mechanisms in the full range of criticality. To do so, a natural way is to incorporate a competition mechanism into the model. This consideration has been inspired by the work of Lambert \cite{Lam05}, who defined a \textit{logistic growth branching process} to model the pairwise competition between individuals in the population. The process was constructed in \cite{Lam05} by a random time change from a spectrally positive Ornstein--Uhlenbeck type process. More general branching population systems with competition were studied in depth by Berestycki et al.\ \cite{BFF18} and Pardoux \cite{Par16}; {\rv see also Foucart \cite{Fou19} and Friesen et al.\ \cite{FJKR23}. Our aim here} is to understand whether and how a strong competition could balance the branching and the immigration to guarantee the exponential ergodicity. {\rrv These are not clear as the process may have infinite first moments, which means the population could be very large. In fact, to study the ergodic behavior we should first understand whether the competition could prevent the population from exploding at a finite time when \eqref{eq1.4} is not satisfied.} The emphasis here is the interplay among the branching, immigration and competition mechanisms.

\subsection{Main results}

Suppose that {\rrv $(\itOmega,\mcr{F},\mcr{F}_t,\mbf{P})$ is a filtered probability space} satisfying the usual hypotheses. Let $\{L(\d s,\d u)\}$ be a spectrally one-sided time-space $(\mcr{F}_t)$-L\'evy white noise with the L\'evy--It\^o decomposition:
 \beqlb\label{eq1.10}
L(\d s,\d u)\ar=\ar W(\d s,\d u) - b\d s\d u + \int_0^1 z \tilde{M}(\d s,\d z,\d u) + \int_1^\infty z M(\d s,\d z,\d u), \nnm
 \eeqlb
where $W(\d s,\d u)$ is a Gaussian white noise on $(0,\infty)^2$ based on $2c\d s\d u$ and $M(\d s,\d z,\d u)$ is a Poisson random measure on $(0,\infty)^3$ with intensity $\d s\mu(\d z)\d u$ and compensated measure $\tilde{M}(\d s,\d z,\d u):= M(\d s,\d z,\d u) - \d s\mu(\d z)\d u$. Here and in the sequel, we make the convention that
 \beqnn
\int_y^x = -\int_x^y = \int_{(y,x]} ~\mbox{and}~ \int_x^\infty = \int_{(x,\infty)}, \quad x\ge y\in \mbb{R}.
 \eeqnn
Let $\{\eta(t): t\ge 0\}$ be an $(\mcr{F}_t)$-subordinator defined by
 \beqlb\label{eq1.11}
\eta(t)= \beta t + \int_0^t\int_0^\infty z N(\d s,\d z),
 \eeqlb
where $N(\d s,\d z)$ is a Poisson random measure on $(0,\infty)^2$ with intensity $\d s\nu(\d z)$. Suppose that $\{L(\d s,\d u)\}$ and $\{\eta(t)\}$ are independent of each other. Let $g$ be a \textit{competition mechanism}, which by definition is a nondecreasing and continuous function on $[0,\infty)$ satisfying $g(0)= 0$. For any $\mcr{F}_0$-measurable random variable $x(0)\ge 0$, we consider the stochastic equation:
 \beqlb\label{eq1.12}
x(t) \ar=\ar x(0) + \int_0^t\int_0^{x(s-)} L(\d s,\d u) - \int_0^t g(x(s)) \d s + \eta(t).
 \eeqlb
We say an $(\mcr{F}_t)$-adapted c\`{a}dl\`{a}g process $\{x(t): t\ge 0\}$ taking values in $\bar{\mbb{R}}_+:= [0,\infty]$ is a \textit{solution} to \eqref{eq1.12} if the equation holds almost surely when $t$ is replaced by $t\land \zeta_n$ for each $t\ge 0$ and $n\ge 1$, where $\zeta_n= \inf\{t\ge 0: x(t)\ge n\}$ with $\inf\emptyset = \infty$. We call $\zeta:= \lim_{n\to \infty} \zeta_n$ the \textit{lifetime} of $\{x(t)\}$ and make the convention that $x(t)= \infty$ for $t\ge \zeta$. If $\zeta= \infty$ almost surely, we say the process $\{x(t)\}$ is \textit{conservative}. We shall prove that there is a pathwise unique solution to \eqref{eq1.12}; see Theorem~\ref{th2.1}. Then the solution $\{x(t)\}$ is a strong Markov process in $\bar{\mbb{R}}_+$. Let $(P_t)_{t\ge 0}$ be the transition semigroup of $\{x(t)\}$. We shall call any Markov process with transition semigroup $(P_t)_{t\ge 0}$ a \textit{continuous-state branching process with immigration and competition} (CBIC-process). Constructions of CB- and CBI-processes in terms of similar stochastic equations were suggested by Bertoin and Le~Gall \cite{BeL06} and Dawson and Li \cite{DaL06, DaL12}. The CBIC-process extends the population models of Berestycki et al.\ \cite{BFF18} and Pardoux \cite{Par16} by the additional immigration structure. A more general continuous-state population model in random environments was introduced by Palau and Pardo \cite{PaP18} by solving a stochastic equation driven by Brownian motions and Poisson random measures.

We need to introduce some concepts in order to present our main results. Let $C^2(\mbb{R}_+)$ be the linear space of twice continuously differentiable functions on $\mbb{R}_+$. {\rv For $x\ge 0$ and} $f\in C^2(\mbb{R}_+)$ write
 \beqlb\label{eq1.13}
Lf(x) \ar=\ar cxf^{\prime\prime}(x) + x\int_0^\infty \big[\itDelta_zf(x) - zf^\prime(x)\I_{\{z\le 1\}}\big] \mu(\d z) \\
 \ar\ar
+\, [\beta-bx-g(x)]f^\prime(x) + \int_0^\infty \itDelta_zf(x) \nu(\d z), \nnm
 \eeqlb
where
 \beqnn
\itDelta_zf(x)= f(x+z) - f(x).
 \eeqnn
Let $\mcr{D}(L)$ denote the linear space consisting of functions $f\in C^2(\mbb{R}_+)$ such that the two integrals on the right-hand side of \eqref{eq1.13} are convergent and define continuous functions on $\mbb{R}_+$. Let $C_b^2(\mbb{R}_+)$ be the space of bounded and continuous functions on $\mbb{R}_+$ with bounded and continuous derivatives up to the second order. Then $C_b^2(\mbb{R}_+)\subset \mcr{D}(L)$. In general, we allow $f$ and $Lf$ to be unbounded functions. We shall see that $(L,\mcr{D}(L))$ is a restriction of the \textit{generator} of the CBIC-process.

Given a nonnegative Borel function $V$ on $\mbb{R}_+$, we denote by $\mcr{P}_V(\mbb{R}_+)$ the set of Borel probability measures $\gamma$ on $\mbb{R}_+$ such that
 \beqnn
\int_{\mbb{R}_+}V(x)\gamma(\d x)< \infty.
 \eeqnn
Let $W_V$ be the \textit{$V$-weighted total variation distance} on $\mcr{P}_V(\mbb{R}_+)$ defined by
 \beqlb\label{eq1.14}
W_V(\gamma,\eta)= \int_{\mbb{R}_+}[1+V(x)]|\gamma-\eta|(\d x),
 \quad
\gamma,\eta\in \mcr{P}_V(\mbb{R}_+),
 \eeqlb
where $|\cdot|$ denotes the total variation measure. We shall see that $W_V$ is actually the \textit{Wasserstein distance} determined by the metric
 \beqlb\label{eq1.15}
d_V(x,y)= [2+V(x)+V(y)]\I_{\{x\neq y\}}, \quad x,y\in \mbb{R}_+.
 \eeqlb
More precisely, we have
 \beqlb\label{eq1.16}
W_V(\gamma,\eta)
 =
\inf_{\pi\in \mcr{C}(\gamma,\eta)}\int_{\mbb{R}_+^2} d_V(x,y) \pi(\d x,\d y),
 \eeqlb
where $\mcr{C}(\gamma,\eta)$ is the collection of all probability measures on $\mbb{R}_+^2$ with marginals $\gamma$ and $\eta$; see Lemma~\ref{th5.3}. The consideration of this distance was inspired by Hairer and Mattingly \cite{HaM11}; see also \cite{EGZ19, LMW21}. In particular, if $V\equiv 0$, then $W_V$ reduces to the total variation distance. The other two frequently used weight functions are given by
 \beqlb\label{eq1.17}
V_1(x)= x ~~\mbox{and}~~ V_{\log}(x)= \log(1+x), \quad x\ge 0.
 \eeqlb

We say {\rv a conservative CBIC-process} $\{x(t): t\ge 0\}$ or its transition semigroup $(P_t)_{t\ge 0}$ is \textit{exponentially ergodic} in the distance $W_V$ with \textit{rate} $\lambda_*>0$ if it possesses a unique stationary probability distribution $\gamma$ and there is a nonnegative function $\eta\mapsto C(\eta)$ on $\mcr{P}_V(\mbb{R}_+)$ such that
 \beqlb\label{eq1.18}
W_V(\gamma,\eta P_t)\le C(\eta)\e^{-\lambda_*t},
 \quad
t\ge 0, \eta\in \mcr{P}_V(\mbb{R}_+).
 \eeqlb
{\rv The exponential ergodicity} \eqref{eq1.18} follows by standard arguments if there is a constant $K\ge 0$ such that
 \beqlb\label{eq1.19}
W_V(P_t(x,\cdot),P_t(y,\cdot))\le K\e^{-\lambda_*t}d_V(x,y), \quad t\ge 0.
 \eeqlb

{\rv Given $\sigma$-finite measures $\mu$ and $\nu$ on $\mbb{R}$, we write $\mu\land \nu$ for $\mu - (\mu-\nu)_+= \nu - (\nu-\mu)_+$, where the subscript ``$+$'' stands for the upper variation of the signed measure in its Jordan decomposition. Let $\mu*\nu$ denote the convolution of $\mu$ and $\nu$ defined by, for all positive Borel functions $f$ on $\mbb{R}$,
 \beqnn
\int_{\mbb{R}} f(z) (\mu*\nu)(\d z)
 =
\int_{\mbb{R}}\mu(\d x)\int_{\mbb{R}} f(x+y)\nu(\d x).
 \eeqnn
For a nonnegative function $V\in \mcr{D}(L)$ and constants $C_0,C_1> 0$ let us consider the inequality}
 \beqlb\label{eq1.20}\rv
LV(x)\le C_0 - C_1V(x), \quad x\ge 0.
 \eeqlb

\bcondition\label{cd1.1} There exists a constant $\lambda_0>0$ such that $\itPsi(\lambda_0)> 0$ and $\itPhi(\lambda_0)> 0$. \econdition

\bcondition\label{cd1.2} {\rv One of the following conditions {\rrv is satisfied:} (i)~Grey's condition \eqref{eq1.9}; (ii)~for some constants $c_0> 0$ and $\kappa_0> 0$,}
 \beqlb\label{eq1.21}\rv
\kappa(x):= [\mu\land(\delta_x*\mu)](0,\infty) + [\nu\land(\delta_x*\nu)](0,\infty)\ge \kappa_0, \quad |x|\le c_0.
 \eeqlb
\econdition

\bcondition\label{cd1.3} There is a nonnegative {\rv function $V\in \mcr{D}(L)$ satisfying \eqref{eq1.20} and} $V(x)\to \infty$ as $x\to \infty$. \econdition

{\rrv We now present the main result of this paper.}

\btheorem\label{th1.1} Suppose that Conditions~\ref{cd1.1}, \ref{cd1.2} and~\ref{cd1.3} are satisfied. Then the CBIC-process is conservative and there are constants $K> 0$ and $\lambda_*> 0$ such that \eqref{eq1.19} holds. Consequently, the CBIC-process is exponentially ergodic in the $V$-weighted total variation distance {\rv with rate $\lambda_*>0$}. \etheorem

The advantage of Theorem~\ref{th1.1} is that it works for general branching mechanisms without criticality restriction. In particular, it applies to all the stable branching mechanisms given by \eqref{eq1.6}. Furthermore, the proof of the theorem actually provides a way of {\rv finding the exponential ergodicity rate} $\lambda_*> 0$, which is important in applications; see Remark~\ref{about_lambda}. {\rv Here, Condition~\ref{cd1.1} is introduced} to avoid some extreme cases of \eqref{eq1.12}. If the condition fails, then either subordinator $\eta(t)$ vanishes or the L\'evy field $L(\d s,\d u)$ only has nonnegative increments. In any of those cases, the immigration essentially plays no role in the ergodic behavior. {\rv Conditions like \eqref{eq1.21}} have been considered in the study of ergodicities of Ornstein-Uhlenbeck type processes with nonlinear drift; see, e.g., \cite{LiW20, LMW21, LuW19, Maj17, ScW12}. Our condition \eqref{eq1.21} is actually weaker than the corresponding assumptions for exponential ergodicities in those previous papers, where one usually required $\kappa(x)\to \infty$ as $x\to 0$. In particular, the condition is satisfied if $\mu(\d z)$ or $\nu(\d z)$ is bounded below by the measure $r\I_{(u,v)}(z)\d z$ for some $r>0$ and $v> u\ge 0$. The existence of {\rv a function $V\in \mcr{D}(L)$ with the property like \eqref{eq1.20} has become a standard assumption in the study of {\rrv uniqueness and ergodicity problems} of Markov processes following Chen \cite{Che86a, Che86b, Che91, Che04}; see also Down et al.\ \cite{DMT95} and Meyn and Tweedie \cite{MeT92, MeT93a, MeT93b}. For the CBIC-process, the condition means roughly that} the competition mechanism wins its confrontation against the branching and immigration when the population is large. Typically, this is equivalent to a simple growth condition for the competition mechanism.

\bproposition\label{th1.2} Suppose that the measures $\mu$ and $\nu$ satisfy the integrability condition
 \beqlb\label{eq1.22}
\int_1^\infty z \mu(\d z) + \int_1^\infty z \nu(\d z)< \infty.
 \eeqlb
Then {\rv $V_1\in \mcr{D}(L)$ satisfies \eqref{eq1.20}} if and only if
 \beqlb\label{eq1.23}
\liminf_{x\to \infty} \frac{g(x)}{x} + b - \int_1^\infty z \mu(\d z)> 0.
 \eeqlb
\eproposition

\bproposition\label{th1.3} Suppose that the measures $\mu$ and $\nu$ satisfy the integrability condition
 \beqlb\label{eq1.24}
\int_1^\infty \log(1+z) \mu(\d z) + \int_1^\infty \log(1+z) \nu(\d z)< \infty.
 \eeqlb
Then {\rv $V_{\log}\in \mcr{D}(L)$ satisfies \eqref{eq1.20}} if and only if
 \beqlb\label{eq1.25}
\liminf_{x\to \infty}\bigg\{\frac{g(x)}{x\log x} - \frac{x}{\log x}\int_1^\infty \log\Big(1+\frac{z}{1+x}\Big) \mu(\d z)\bigg\}> 0.
 \eeqlb
\eproposition

Under the integrability condition \eqref{eq1.22}, the CBIC-process {\rrv has finite first moments.} In this case, condition \eqref{eq1.23} means the competition mechanism should grow at least linearly with a sufficiently large rate. In particular, if $g(x)= ax$ for $x\ge 0$, the CBIC-process reduces to a CBI-process with branching mechanism $\lambda\mapsto a\lambda + \itPsi(\lambda)$ and \eqref{eq1.23} simply means the branching mechanism should be subcritical. This is consistent with the subcriticality condition in \cite[Theorem~2.5]{LiM15}, but here we do not need to assume Grey's condition. The conditions of Proposition~\ref{th1.3} apply to supercritical branching mechanisms including those with infinite first moments but finite logarithmic moments.

\subsection{Stable branching CBIC-processes}

These processes are solutions of stochastic differential equations driven by spectrally one-sided L\'evy processes. Let $m_\alpha$ be the $\sigma$-finite measure on $(0,\infty)$ defined by
 \beqlb\label{eq1.26}
m_\alpha(\d z)= \left\{\begin{array}{ll}
 (\alpha-1)\itGamma(2-\alpha)^{-1} z^{-1-\alpha}\d z, & 1< \alpha< 2, \cr
 z^{-2}\d z, & \alpha= 1, \cr
 \itGamma(1-\alpha)^{-1}z^{-1-\alpha}\d z, & 0< \alpha< 1,
\end{array}\right.
 \eeqlb
where $\itGamma$ is the Gamma function. Let $\{M_\alpha(\d s,\d z)\}$ be a time-space $(\mcr{F}_t)$-Poisson random measure on $(0,\infty)^2$ with intensity $\d sm_\alpha(\d z)$ and compensated measure $\tilde{M}_\alpha(\d s,\d z)$. Let $\{z_\alpha(t)\}$ be the spectrally positive $\alpha$-stable L\'evy process defined by
 \beqlb\label{eq1.27}
z_\alpha(t)= \left\{\begin{array}{ll}
 \displaystyle\int_0^t\int_0^\infty z \tilde{M}_\alpha(\d s,\d z), & 1< \alpha< 2, \smallskip\cr
 \displaystyle\int_0^t\int_0^1 z \tilde{M}_\alpha(\d s,\d z) + \int_0^t\int_1^\infty z M_\alpha(\d s,\d z), \quad & \alpha= 1, \smallskip\cr
 \displaystyle\int_0^t\int_0^\infty z M_\alpha(\d s,\d z), & 0< \alpha< 1.
\end{array}\right.
 \eeqlb
Let $\{B(t)\}$ be a standard $(\mcr{F}_t)$-Brownian motion and let $\{\eta(t)\}$ be the $(\mcr{F}_t)$-subordinator defined by \eqref{eq1.11}. Suppose that $\{B(t)\}$, $\{z_\alpha(t)\}$ and $\{\eta(t)\}$ are independent of each other. Given any $\mcr{F}_0$-measurable random variable $x(0)\ge 0$, we shall construct the CBIC-process $\{x(t): t\ge 0\}$ with stable branching mechanism defined by \eqref{eq1.6} in terms of the stochastic differential equations, for $\alpha= 1$,
 \beqlb\label{eq1.28}
\d x(t)= \sqrt{2cx(t)}\d B(t) + \sigma x(t-)\d z_1(t) + \sigma x(t)\log(\sigma x(t))\d t - ax(t)\d t - g(x(t))\d t + \d \eta(t)
 \eeqlb
and, for $\alpha\neq 1$,
 \beqlb\label{eq1.29}
\quad \d x(t) = \sqrt{2cx(t)}\d B(t) + \sqrt[\alpha]{\alpha\sigma x(t-)}\d z_\alpha(t) - ax(t)\d t - g(x(t))\d t + \d \eta(t);
 \eeqlb
see Theorem~\ref{th2.4}. By applying Theorem~\ref{th1.1} and Propositions~\ref{th1.2} and~\ref{th1.3} to the stable branching mechanisms, we obtain the following:

\bcorollary\label{th1.4} Suppose that the branching mechanism is given by \eqref{eq1.6} with $1< \alpha< 2$. In addition, assume that
 \beqnn
\int_1^\infty z \nu(\d z)< \infty
 ~~\mbox{and}~~
\liminf_{x\to \infty} \frac{g(x)}{x}> -a.
 \eeqnn
Then the CBIC-process is exponentially ergodic in the $V_1$-weighted total variation distance. \ecorollary

\bcorollary\label{th1.5} Suppose that the branching mechanism is given by \eqref{eq1.6} with $0< \alpha\le 1$ and that $a>0$ or $c>0$ for $0< \alpha< 1$. In addition, assume that
 \beqnn
\int_1^\infty \log(1+z) \nu(\d z)< \infty
 \eeqnn
and
 \beqnn
\left\{\begin{array}{ll}
 \displaystyle\liminf_{x\to \infty} \frac{g(x)}{x^{2-\alpha}}> \frac{\sigma\pi}{\itGamma(1-\alpha)\sin(\alpha\pi)}, &\quad 0< \alpha< 1, \cr
 \displaystyle\liminf_{x\to \infty} \frac{g(x)}{x\log x}> \sigma, &\quad \alpha= 1.
\end{array}\right.
 \eeqnn
Then the CBIC-process is exponentially ergodic in the $V_{\log}$-weighted total variation distance.
\ecorollary

In the situation of Corollary~\ref{th1.5}, the CBIC-process has {\rrv infinite first moments,} so the result has no counterpart in the setting of CBI-processes, where exponential ergodicities can only be considered for subcritical branching mechanisms. The assumption that $a>0$ or $c>0$ for $0< \alpha< 1$ in the corollary comes from Condition~\ref{cd1.1}. If this were not satisfied, there would not be sufficient fluctuations in the CBIC-process as the two remaining noise terms in \eqref{eq1.29} are both nondecreasing. As we explained before, if $\sigma>0$ and $0<\alpha< 1$, then condition \eqref{eq1.4} fails and the corresponding CBI-process without competition explodes at a finite time.

\subsection{{\rv More examples}}

The following examples concern further applications of the results presented above and show the conditions we introduce are sharp in typical situations.

{\rv

\bexample\label{ex1.0a} Consider the CBI-process with branching and immigration mechanisms given respectively by \eqref{eq1.1} and \eqref{eq1.7} with $b> 0$, $\beta> 0$ and $c=\nu(0,\infty)= 0$. We assume that $\mu(\d z)$ is given by
 \beqnn
\mu(\d z)= \alpha \sum^\infty_{k=1} m_\alpha [k, k+1)\delta_{k+1}(\d z) + \alpha \sum_{k=1}^\infty m_\alpha[1/(k+1),1/k)\delta_{1/k}(\d z),
 \eeqnn
where $1< \alpha< 2$ and $m_\alpha(\d z)$ is the $\sigma$-finite measure on $(0,\infty)$ defined by \eqref{eq1.26}. Then the branching mechanism satisfies Grey's condition \eqref{eq1.9} since
 \beqnn
\itPsi(\lambda)> b\lambda + \lambda^\alpha= b\lambda + \alpha \int_0^\infty \big(\e^{-\lambda z}-1+\lambda z\big)m_\alpha(\d z).
 \eeqnn
By Proposition~\ref{th1.2}, the CBI-process is exponentially ergodic relative to the $V_1$-weighted total variation distance. {\rrv The exponential ergodicity} in the total variation distance of the process can also {\rrv be derived from} Li and Ma \cite[Theorem~2.5]{LiM15}. In this case, Condition~\ref{cd1.2}-(ii) is not satisfied because $[\mu\land(\delta_x*\mu)](0,\infty)= 0$ when $x$ is any irrational number. \qed \eexample

}

\bexample\label{ex1.1} Consider the CBI-process with branching and immigration mechanisms given respectively by \eqref{eq1.1} and \eqref{eq1.7} with $b> 0$, $\beta> 0$ and $c=\nu(0,\infty)= 0$. In addition, assume that $\mu(\d z)= r\I_{(u,v)}(z)\d z$ for some $r>0$ and $0\le u< v\le 1$. {\rv In this case, the branching mechanism satisfies Condition~\ref{cd1.2}-(ii).} By Proposition~\ref{th1.2}, the CBI-process is exponentially ergodic relative to the $V_1$-weighted total variation distance. {\rv However, the branching mechanism} does not satisfy Grey's condition \eqref{eq1.9} as
 \beqnn
\lim_{\lambda\to \infty} \lambda^{-1}\itPsi(\lambda)
 =
b + \int_0^1 z \mu(\d z)= b + \frac{r}{2}(v^2-u^2).
 \eeqnn
Therefore the exponentially ergodicity in the total variation distance of the process does not follow from Li and Ma \cite[Theorem~2.5]{LiM15}. One can also see that the conditions of Li and Wang \cite[Theorem~1.1]{LiW20} are not satisfied. \qed \eexample

\bexample\label{ex1.2} Consider the CBI-process with branching and immigration mechanisms given by \eqref{eq1.1} and \eqref{eq1.7}, respectively. Suppose that \eqref{eq1.22} holds and the branching mechanism is critical, that is,
 \beqnn
b_0:= b - \int_1^\infty z \mu(\d z)= 0.
 \eeqnn
It is known that
 \beqnn
\int_{\mbb{R}_+} y P_t(0,\d y)= t\beta+ t\int_0^\infty z\nu(\d z);
 \eeqnn
see, e.g., Li \cite[p.33]{Li20}. The right-hand side tends to infinity as $t\to \infty$, so the CBI-process is not exponentially ergodic in the $V_1$-weighted total variation distance. This shows the conditions of Corollary~\ref{th1.4} are sharp. \qed \eexample

\bexample\label{ex1.3} Let $\itPsi$ be the stable branching mechanisms given by \eqref{eq1.6}. By the proof of Theorem~\ref{th2.4}, the CBIC-process has generator $L$ defined by \eqref{eq1.13} with $b= a + \sigma h_\alpha$ and $\mu(\d z)= \alpha\sigma m_\alpha(\d z)$, where $m_\alpha$ and $h_\alpha$ are defined by \eqref{eq1.26} and \eqref{eq2.6}, respectively. Let us consider the case where $\alpha=1$ and $b= c= \nu(0,\infty)= 0$. It is easy to show that
\beqnn
LV_{\log}(x)= \frac{\sigma x}{1+x}\big[1+\log(1+x)\big] - \frac{g(x)}{1+x} + \frac{\beta}{1+x}.
\eeqnn
If $\sigma\beta>0$ and $g(x)= \sigma x\log(1+x)$ for $x\ge 0$, there is a constant $\lambda> 0$ such that $LV_{\log}(x)\ge \lambda$, and so
 \beqnn
\int_{\mbb{R}_+} V_{\log}(y) P_t(x,\d y)
 =
V_{\log}(x) + \int^t_0 P_s LV_{\log}(x)\d s
 \ge
V_{\log}(x) + \lambda t.
 \eeqnn
Consequently, the CBIC-process is not exponentially ergodic in the $V_{\log}$-weighted total variation distance. When $0< \alpha< 1$, a similar example can be given by considering the function $g(x)= [\itGamma(1-\alpha)\sin(\alpha\pi)]^{-1}\sigma\pi x^{2-\alpha}$. Then the conditions of Corollary~\ref{th1.5} are sharp. \qed \eexample

\subsection{The approach of coupling and distance}

The characterizations like \eqref{eq1.3} and \eqref{eq1.8} are not available for the CBIC-process. Our proof of Theorem~\ref{th1.1} is based on the approach of coupling and distance. Those techniques have played important roles in the ergodic theory of Markov processes. We refer to Chen \cite{Che04, Che05} for systematical treatments of the techniques. In particular, a number of results on the exponential ergodicity have been obtained by this method for \textit{Ornstein-Uhlenbeck type} processes with nonlinear drifts defined by stochastic differential equations of the form:
 \beqlb\label{eq1.30}
\d x(t) = \d L(t) - g(x(t))\d t, \quad t\ge 0,
 \eeqlb
where $\{L(t): t\ge 0\}$ is a L\'evy process and the drift coefficient $g$ is sufficiently regular so that there is a unique solution to the equation. When the driving noise is a Brownian motion, the exponential ergodicity and related problems were studied by Eberle and his coauthors \cite{Ebe11, Ebe16, EGZ19} and Luo and Wang \cite{LuW16}. For discontinuous L\'evy noises, the ergodicity problem was investigated by Liang et al.\ \cite{LMW21}, Luo and Wang \cite{LuW19}, Schilling and Wang \cite{ScW12} and Majka \cite{Maj17}. The feature of those processes is that the random noise in \eqref{eq1.30} is both temporarily and spatially homogeneous.

The main difficulty of the proof of Theorem~\ref{th1.1} is that the branching fluctuations become small and rare when the process is close to zero, which is clear from \eqref{eq1.12} or \eqref{eq1.13}. A similar phenomenon has been noticed in nonlinear branching processes by Li and Wang \cite{LiW20}; see also \cite{Li19, LYZ19}. The difficulty was treated in \cite{LiW20} by introducing conditions on the asymptotics at zero of the branching coefficients. We cannot do that in our setting here since the branching coefficients here have to be linear. The key of our proof is the design of a Markov coupling process and a nonsymmetric control function. More precisely, we construct a two-dimensional \textit{Markov coupling process} $\{(X_t,Y_t): t\ge 0\}$ on $\mbb{R}_+^2$, where both $\{X_t: t\ge 0\}$ and $\{Y_t: t\ge 0\}$ are Markov processes with transition semigroup $(P_t)_{t\ge 0}$. The coupling process satisfies $X_{T+t}= Y_{T+t}$ for every $t\ge 0$, where $T= \inf\{t\ge 0: X_t= Y_t\}$ with $\inf\emptyset= \infty$ by convention is the \textit{succeeding} or \textit{coupling time}. We need to make the coupling succeed as early as possible, which is realized by reflecting the noises in suitable ways; see Remark~\ref{th3.6}. Then we define a non-symmetric function $G_0$ with the \textit{exponential contraction property}: for $t\ge 0$ and $(x,y)\in \mbb{R}_+^2$,
 \beqlb\label{eq1.31}
\tilde{P}_tG_0(x,y)\le \e^{-\lambda_*t}G_0(x,y),
 \eeqlb
where $(\tilde{P}_t)_{t\ge 0}$ is the transition semigroup of the coupling process. The function $G_0$ should \textit{control} the distance $d_V$ in the sense that there are constants $c_2\ge c_1> 0$ such that, for $(x,y)\in \mbb{R}_+^2$,
 \beqlb\label{eq1.32}
c_1G_0(x,y)\le d_V(x,y)\le c_2G_0(x,y).
 \eeqlb
Then we deduce the estimate \eqref{eq1.19} from \eqref{eq1.31} and \eqref{eq1.32}. The non-symmetric control function makes it possible to identify the dominating factor from the branching, immigration and competition mechanisms in different parts of the space, which is important in establishing the estimate \eqref{eq1.31}. By those devices we not only overcome the difficulty caused by the degeneracy of the branching coefficients but also weaken the fluctuation conditions introduced for the exponential ergodicity of the Ornstein-Uhlenbeck type processes; see, e.g., \cite{LiW20, LMW21, LuW19, Maj17, ScW12}. In fact, the assumptions on the drift coefficient in Corollaries~\ref{th1.4} and~\ref{th1.5} are also weaker than those in the previous papers, where lower bounds of the increment $g(x)-g(y)$ for $x\ge y\ge0$ were required. We believe that the method used here could be adapted to more general classes of processes with L\'evy driving noises.

The rest of the paper is organized as follows. In Section~2, the existence and uniqueness of solutions to the stochastic equations \eqref{eq1.12}, \eqref{eq1.28} and \eqref{eq1.29} are proved. In Section~3, we give the construction of the Markov coupling process. In Section~4, some necessary estimates for the coupling generator are established. The proofs of the ergodicity results are given in Section~5.

\smallskip

{\rv\noindent\textbf{Acknowledgement.} We would like to express our sincere thanks to an anonymous referee for helpful comments and suggestions, which have led to a number of improvements of the paper.}

\section{Stochastic equations of CBIC-processes}

 \setcounter{equation}{0}

In this section, we prove the existence and pathwise uniqueness of solutions to the stochastic equations \eqref{eq1.12}, \eqref{eq1.28} and \eqref{eq1.29}, which give constructions of the CBIC-processes.

\btheorem\label{th2.1} For any $\mcr{F}_0$-measurable random variable $x(0)\ge 0$, there is a pathwise unique solution $\{x(t): t\ge 0\}$ to \eqref{eq1.12}. \etheorem

\begin{proof} By Theorem~2.5 of Dawson and Li \cite{DaL12}, for each $n\ge 1$ there is a pathwise unique solution $\{x_n(t): t\ge 0\}$ to
 \beqlb\label{eq2.1}
\quad x(t) \ar=\ar x(0) + \int_0^t\int_0^{x(s-)}W(\d s,\d u) + \int_0^t [\beta-bx(s)-g(x(s))] \d s \\
 \ar\ar
+ \int_0^t\int_0^1 \int_0^{x(s-)} z \tilde{M}(\d s,\d z,\d u) + \int_0^t\int_1^\infty \int_0^{x(s-)} (z\land n) M(\d s,\d z,\d u) \nnm\\
 \ar\ar
+ \int_0^t\int_0^\infty (z\land n) N(\d s,\d z). \nnm
 \eeqlb
Let $\zeta_n= \inf\{t\ge 0: x_n(t)\ge n\}$. From the pathwise uniqueness for \eqref{eq2.1} it follows that $\zeta_n$ is nondecreasing in $n\ge 1$ and $x_{n+1}(t)= x_n(t)$ for $0\le t< \zeta_n$. Clearly, the pathwise unique solution $\{x(t): t\ge 0\}$ to \eqref{eq1.12} is given by $x(t)= x_n(t)$ for $0\le t< \zeta_n$ and $x(t)= \infty$ for $t\ge \zeta:= \lim_{n\to \infty}\zeta_n$. \end{proof}

{\rrv By the above theorem we have given a construction of the CBIC-process.} The next result justifies the fact that the operator $(L,\mcr{D}(L))$ defined by \eqref{eq1.13} is \textit{a restriction of the generator} of the process.

\btheorem\label{th2.2} Let $\{x(t): t\ge 0\}$ be the pathwise unique solution to \eqref{eq1.12} and let $\zeta_n= \inf\{t\ge 0: x(t)\ge n\}$. Then for any $n\ge 1$ and $f\in \mcr{D}(L)$ we have
 \beqlb\label{eq2.2}
f(x(t\land \zeta_n))= f(x(0)) + \int_0^{t\land \zeta_n} Lf(x(s)) \d s + M_n(t),
 \eeqlb
where $\{M_n(t): t\ge 0\}$ is a martingale defined by
 \beqlb\label{eq2.3}
\quad M_n(t)\ar=\ar \int_0^{t\land \zeta_n}\int_0^{x(s-)}f'(x(s-))W(\d s,\d u) + \int_0^{t\land \zeta_n}\int_0^\infty \itDelta_zf(x(s-)) \tilde{N}(\d s,\d z) \\
 \ar\ar
+ \int_0^{t\land \zeta_n}\int_0^\infty\int_0^{x(s-)} \itDelta_zf(x(s-)) \tilde{M}(\d s,\d z,\d u). \nnm
 \eeqlb
\etheorem

\begin{proof} By \eqref{eq1.12} the process $\{x(t): t\ge 0\}$ is a semimartingale. For any $f\in \mcr{D}(L)$, we can use It\^o's formula to see that
 \beqnn
f(x(t\land \zeta_n))\ar=\ar f(x(0)) - \int_0^{t\land \zeta_n} f'(x(s))g(x(s)) \d s + \int_0^{t\land \zeta_n}\int_0^{x(s-)} f'(x(s-)) L(\d s,\d u) \cr
 \ar\ar
+ \int_0^{t\land \zeta_n} f'(x(s-)) \eta(\d s) + c\int_0^{t\land \zeta_n} f''(x(s))x(s) \d s \cr
 \ar\ar
+ \int_0^{t\land \zeta_n}\int_0^\infty\int_0^{x(s-)} \big[\itDelta_zf(x(s-)) - zf'(x(s-))\big] M(\d s,\d z,\d u) \cr
 \ar\ar
+ \int_0^{t\land \zeta_n}\int_0^\infty \big[\itDelta_zf(x(s-)) - zf'(x(s-))\big] N(\d s,\d z) \cr
 \ar=\ar
f(x(0)) - \int_0^{t\land \zeta_n} f'(x(s))g(x(s)) \d s + \int_0^{t\land \zeta_n}\int_0^{x(s-)} f'(x(s-)) W(\d s,\d u) \cr
 \ar\ar
+ \int_0^{t\land \zeta_n} \big[\beta-b x(s)\big] f'(x(s)) \d s + c\int_0^{t\land \zeta_n} x(s)f''(x(s)) \d s \cr
 \ar\ar
+ \int_0^{t\land \zeta_n}\int_0^\infty\int_0^{x(s-)} \itDelta_zf(x(s-)) \tilde{M}(\d s,\d z,\d u) \cr
 \ar\ar
+ \int_0^{t\land \zeta_n}x(s)\d s\int_0^\infty \big[\itDelta_zf(x(s)) - zf'(x(s))\I_{\{z\le 1\}}\big] \mu(\d z) \cr
 \ar\ar
+ \int_0^{t\land \zeta_n}\int_0^\infty \itDelta_zf(x(s-)) \tilde{N}(\d s,\d z)
+ \int_0^{t\land \zeta_n}\d s\int_0^\infty \itDelta_zf(x(s)) \nu(\d z).
 \eeqnn
Then \eqref{eq2.2} holds with $M_n(t)$ defined by \eqref{eq2.3}. Since $x(s-)\le n$ for $0< s\le \zeta_n$, it is easy to show that $\{M_n(t): t\ge 0\}$ is a martingale. \end{proof}

\bproposition\label{th2.3} Suppose that Condition~\ref{cd1.3} is satisfied. Then the solution $\{x(t): t\ge 0\}$ to \eqref{eq1.12} is conservative and
 \beqlb\label{eq2.4}
\mbf{E}[V(x(t))]\le \mbf{E}[V(x(0))]\e^{-C_1t} + C_0C_1^{-1}\big(1-\e^{-C_1t}\big), \quad t\ge 0.
 \eeqlb
\eproposition

\begin{proof} There is no loss of generality to assume $\mbf{E}[V(x(0))]< \infty$. {\rv Recall that $\zeta_n= \inf\{t\ge 0: x(t)\ge n\}$ and $\zeta= \lim_{n\to \infty}\zeta_n$.} By applying \eqref{eq2.2} to the function $V$ and using integration by parts we have
 \beqlb\label{eq2.5}\rv
\e^{C_1(t\land \zeta_n)}V(x(t\land \zeta_n))
 =
V(x(0)) + \int_0^{t\land \zeta_n} \e^{C_1s}\big[C_1V(x(s))+LV(x(s))\big] \d s + M(t),
 \eeqlb
where
 \beqnn
M(t)\ar=\ar \rv\int_0^{t\land \zeta_n}\int_0^{x(s)}\e^{C_1s}f'(x(s))W(\d s,\d u) + \int_0^{t\land \zeta_n}\int_0^\infty \e^{C_1s}\itDelta_zf(x(s-)) \tilde{N}(\d s,\d z) \cr
 \ar\ar\rv
+ \int_0^{t\land \zeta_n}\int_0^\infty\int_0^{x(s-)} \e^{C_1s}\itDelta_zf(x(s-)) \tilde{M}(\d s,\d z,\d u).
 \eeqnn
{\rv In view of \eqref{eq1.20}, we can take the expectations in both sides of \eqref{eq2.5}} to obtain
 \beqnn
\mbf{E}[\e^{C_1(t\land \zeta_n)}V(x(t\land \zeta_n))]\le \mbf{E}[V(x(0))] + C_0C_1^{-1}\big(\e^{C_1t}-1\big).
 \eeqnn
By Fatou's lemma it follows that
 \beqnn
\mbf{E}[\e^{C_1(t\land \zeta)}V(x(t\land \zeta))]\le \mbf{E}[V(x(0))] + C_0C_1^{-1}\big(\e^{C_1t}-1\big)
 \eeqnn
with the convention $V(\infty)= \infty$. Then $\mbf{P}(\zeta\le t)= 0$ and \eqref{eq2.4} follows. Since $t\ge 0$ was arbitrary, we have $\zeta= \infty$ almost surely. \end{proof}

\btheorem\label{th2.4} For any $\mcr{F}_0$-measurable random variable $x(0)\ge 0$, there are pathwise unique solutions to \eqref{eq1.28} and \eqref{eq1.29}, and the solutions are CBIC-processes with competition mechanism $g$, stable branching mechanisms $\itPsi$ given by \eqref{eq1.6} and immigration mechanism $\itPhi$ given by \eqref{eq1.7}.
\etheorem

\begin{proof} Let $\{x(t)\}$ be the solution to \eqref{eq1.12} with $\{M(\d s,\d z,\d u)\}$ being a Poisson random measure on $(0,\infty)^3$ with intensity $\alpha\sigma\d sm_\alpha(\d z)\d u$, where $m_\alpha(\d z)$ is given by \eqref{eq1.26}. Moreover, we take $b= a + \sigma h_\alpha$, where
 \beqlb\label{eq2.6}
h_\alpha= \left\{\begin{array}{ll}
 -\alpha\itGamma(2-\alpha)^{-1}, &\quad 1< \alpha< 2, \cr
 0, &\quad \alpha= 1, \cr
 \alpha(1-\alpha)^{-1}\itGamma(1-\alpha)^{-1}, &\quad 0< \alpha< 1.
\end{array}\right.
 \eeqlb
If $c>0$, one can see that
 \beqnn
B(t)= \int_0^t\int_0^{x(s)} \frac{1}{\sqrt{2cx(s)}}\I_{\{x(s)> 0\}} W(\d s,\d u) + \frac{1}{\sqrt{2c}} \int_0^t\int_0^1 \I_{\{x(s)= 0\}} W(\d s,\d u)
 \eeqnn
defines a standard Brownian motion and
 \beqnn
\int_0^t\sqrt{2cx(s)}\d B(s)
 =
\int_0^t\int_0^{x(s)} W(\d s,\d u).
 \eeqnn
In the sequel, we assume $\sigma> 0$, for otherwise the proof is simpler. Let $\{M_\alpha(\d s,\d z)\}$ be the random measure on $(0,\infty)^2$ defined by, for $B\in \mcr{B}(0,\infty)$,
 \beqnn
M_\alpha((0,t]\times B)
 \ar=\ar
\int_0^t\int_0^\infty\int_0^{x(s-)} \I_{\{x(s-)>0\}} \I_B\Big(\frac{z}{\sqrt[\alpha]{\alpha\sigma x(s-)}}\Big) M(\d s,\d z,\d u) \cr
 \ar\ar
+ \int_0^t\int_0^\infty\int_0^{1/(\alpha\sigma)} \I_{\{x(s-)=0\}}\I_B(z) M(\d s,\d z,\d u).
 \eeqnn
One can show {\rv as in Li \cite[pp.287--288]{Li22} that} $\{M_\alpha(\d s,\d z)\}$ is a Poisson random measure with intensity $\d sm_\alpha(\d z)$. Let $\{z_\alpha(t)\}$ be the $\alpha$-stable L\'evy processes defined by \eqref{eq1.27}. It is easy to see that
 \beqnn
\int_0^t\sqrt[\alpha]{\alpha\sigma x(s-)}\d z_\alpha(s)
 \ar=\ar
\int_0^t\int_0^1 \int_0^{x(s-)} z \tilde{M}(\d s,\d z,\d u) + \sigma h_\alpha\int_0^t x(s)\d s\cr
 \ar\ar
+ \int_0^t\int_1^\infty \int_0^{x(s-)} z M(\d s,\d z,\d u).
 \eeqnn
Then \eqref{eq1.28} and \eqref{eq1.29} hold. By Fu and Li \cite[Theorem~5.1]{FuL10}, for any $n\ge 1$ there is a pathwise unique solution to the stochastic equation:
 \beqnn
x(t) \ar=\ar x(0) + \int_0^t \sqrt{2cx(s)} \d B(s) + \int_0^t\int_0^1 \sigma x(s-)z \tilde{M}_1(\d s,\d z) \cr
 \ar\ar
+ \int_0^t\int_1^\infty \sigma x(s-)(z\land n) M_1(\d s,\d z) + \sigma\int_0^t x(s)\log(\sigma x(s))\d s \cr
 \ar\ar
- \int_0^t [ax(s)+g(x(s))]\d s + \eta(t).
 \eeqnn
Since $n\ge 1$ is arbitrary, the pathwise uniqueness also holds for \eqref{eq1.28}. A similar argument gives the pathwise uniqueness of \eqref{eq1.29}. These give the results of the theorem; see, e.g., Situ \cite[p.76 and p.104]{Sit05}. \end{proof}

\section{Construction of the coupling process}

 \setcounter{equation}{0}

{\rv Throughout this section, we assume the CBIC-process with an arbitrary initial value is conservative. By Proposition~\ref{th2.3}, this is true if Condition~\ref{cd1.3} is satisfied. We shall give the construction of a Markov coupling of the CBIC-process.} Given $x\in \mbb{R}$ and a general $\sigma$-finite measure $m$ on $(0,\infty)$, we define another $\sigma$-finite measure on $(0,\infty)$ by
 \beqlb\label{eq3.1}
m_x(\d z)= \frac{1}{2}\big[m\land(\delta_x*m)\big](\d z).
 \eeqlb
It is easy to see that
 \beqlb\label{eq3.2}
m_x= \delta_x*m_{-x}
 ~~\mbox{and}~~
m_x(0,0\vee x]= 0.
 \eeqlb
Moreover, we have
 \beqlb\label{eq3.3}
m_x(0,\infty)= m_{-x}(0,\infty)\le \frac{1}{2}m(|x|,\infty).
 \eeqlb
For any $x,y\ge 0$ write $m_{x,y}= m_{y-x} + m_{x-y}$. One can show that both $x\mapsto m_x(\d z)$ and $(x,y)\mapsto m_{x,y}(\d z)$ are kernels. For notational convenience, we also write those kernels as $m(x,\d z)$ and $m(x,y,\d z)$, respectively. We shall use those notations and facts to the L\'evy measures $\mu$ and $\nu$ of the CBIC-process. For $x\in \mbb{R}$ and $z>0$ we set
 \beqlb\label{eq3.4}
\rho_1(x,z)= \frac{\mu(x,\d z)}{\mu(\d z)},
 ~~
r_1(x,z)= \frac{\nu(x,\d z)}{\nu(\d z)}.
 \eeqlb
By \eqref{eq3.1}, the above Radon-Nikodym derivatives exist and satisfy
 \beqnn
\sup_{x\in \mbb{R}}\sup_{z>0}\rho_1(x,z)\le \frac{1}{2},
 \quad
\sup_{x\in \mbb{R}}\sup_{z>0}r_1(x,z)\le \frac{1}{2}.
 \eeqnn
In view of \eqref{eq3.2}, we have
 \beqlb\label{eq3.5}
\rho_1(x,z)= r_1(x,z)= 0, \quad 0<z\le 0\vee x.
 \eeqlb
For $x,y\ge 0$ and $z>0$ let
 \beqnn
\rho_2(x,y,z)= \rho_1(x-y,z) + \rho_1(y-x,z)
 \eeqnn
and
 \beqnn
r_2(x,y,z)= r_1(x-y,z) + r_1(y-x,z).
 \eeqnn

Suppose that $(\itOmega,\mcr{F},\mcr{F}_t,\mbf{P})$ is a filtered probability space satisfying the usual hypotheses. Let $\{W_0(\d s,\d u)\}$ {\rv be a time-space} $(\mcr{F}_t)$-Gaussian white noise on $(0,\infty)^2$ based on $2c\d s\d u$. Let $\{M_0(\d s,\d z,\d u,\d v)\}$ be a time-space $(\mcr{F}_t)$-Poisson random measure on $(0,\infty)^3\times (0,1]$ with intensity $\d s\mu(\d z)\d u\d v$ and compensated measure $\{\tilde{M}_0(\d s,\d z,\d u,\d v)\}$. Let $\{N_0(\d s,\d z,\d v)\}$ be a time-space $(\mcr{F}_t)$-Poisson random measure on $(0,\infty)^2\times (0,1]$ with intensity $\d s\nu(\d z)\d v$. We assume those random noises are independent of each other. Here the extra component $(0,1]$ of the Poisson random measures is introduced for the convenience of disassembling the jumps in the construction of the coupling process. Let $\{\eta_t\}$ be the $(\mcr{F}_t)$-subordinator defined by
 \beqnn
\eta_t= \beta t + \int_0^t\int_0^\infty\int_0^1 z N_0(\d s,\d z,\d v).
 \eeqnn
Let $\{L_0(\d s,\d u)\}$ and $\{L_0^*(\d s,\d u)\}$ be two accompanying spectrally one-sided time-space $(\mcr{F}_t)$-L\'evy white noises defined by
 \beqlb\label{eq3.6}
\quad L_0(\d s,\d u)\ar=\ar W_0(\d s,\d u) - b\d s\d u + \int_{\{0< z\le 1\}}\int_{\{0< v\le 1\}} z \tilde{M}_0(\d s,\d z,\d u,\d v) \\
 \ar\ar
+ \int_{\{1< z< \infty\}}\int_{\{0< v\le 1\}} z M_0(\d s,\d z,\d u,\d v) \nnm
 \eeqlb
and
 \beqlb\label{eq3.7}
\quad L_0^*(\d s,\d u)\ar=\ar -\,W_0(\d s,\d u) - b\d s\d u + \int_{\{0< z\le 1\}}\int_{\{0< v\le 1\}} z \tilde{M}_0(\d s,\d z,\d u,\d v) \\
 \ar\ar
+ \int_{\{1< z< \infty\}}\int_{\{0< v\le 1\}} z M_0(\d s,\d z,\d u,\d v). \nnm
 \eeqlb
By Theorem~\ref{th2.1}, for any $\mcr{F}_0$-measurable random variable $X_0\ge 0$ we can construct a CBIC-process $\{X_t: t\ge 0\}$ by the pathwise unique solution to the stochastic integral equation:
 \beqlb\label{eq3.8}
X_t= X_0 + \int_0^t\int_0^{X_{s-}}L_0(\d s,\d u) - \int_0^t g(X_s)\d s + \eta_t.
 \eeqlb

\btheorem\label{th3.1} {\rv Let $Y_0$ be} an $\mcr{F}_0$-measurable random variable satisfying $X_0\ge Y_0\ge 0$. Then there is a pathwise unique solution $\{(Y_t,\xi_t): t\ge 0\}$ to the system of stochastic equations:
 \beqlb\label{eq3.9}
Y_t\ar=\ar Y_0 + \int_{t\land T}^t\int_0^{Y_{s-}}L_0(\d s,\d u) - \int_0^t g(Y_s)\d s + \eta_t \\
 \ar\ar
+ \int_0^{t\land T}\int_0^{Y_{s-}}L_0^*(\d s,\d u) + \xi_t \nnm
 \eeqlb
and
 \beqlb\label{eq3.10}
\xi_t\ar=\ar \int_0^{t\land T}\int_0^\infty\int_0^{Y_{s-}}\int_0^{\rho_1(Y_{s-}-X_{s-},z)} (X_{s-}-Y_{s-}) M_0(\d s,\d z,\d u,\d v) \\
 \ar\ar
+ \int_0^{t\land T}\int_0^\infty\int_0^{r_1(Y_{s-}-X_{s-},z)} (X_{s-}-Y_{s-}) N_0(\d s,\d z,\d v) \nnm\\
 \ar\ar
+ \int_0^{t\land T}\int_0^\infty\int_0^{Y_{s-}}\int_{\rho_1(Y_{s-}-X_{s-},z)}^{\rho_2(X_{s-},Y_{s-},z)} (Y_{s-}-X_{s-}) M_0(\d s,\d z,\d u,\d v) \nnm\\
\ar\ar
+ \int_0^{t\land T}\int_0^\infty\int_{r_1(Y_{s-}-X_{s-},z)}^{r_2(X_{s-},Y_{s-},z)} (Y_{s-}-X_{s-}) N_0(\d s,\d z,\d v), \nnm
 \eeqlb
where
 \beqlb\label{eq3.11}
T= \inf\{t\ge 0: X_t\le Y_t\}= \inf\{t\ge 0: X_t= Y_t ~\mbox{or}~ X_{t-}= Y_{t-}\}.
 \eeqlb
Moreover, we have $X_{T+t}= Y_{T+t}$ for every $t\ge 0$ if $T< \infty$.
\etheorem

In the proof of the above theorem, we shall see that the pure jump process $\{\xi_t: t\ge 0\}$ has at most a finite number of jumps. In order to give the proof, we make some preparations. By Theorem~\ref{th2.1}, there is a CBIC-process $\{Z_0(t): t\ge 0\}$ defined by the pathwise unique solution to
 \beqlb\label{eq3.12}
Z_0(t)= Y_0 + \int_0^t\int_0^{Z_0(s-)}L_0^*(\d s,\d u) - \int_0^t g(Z_0(s))\d s + \eta_t.
 \eeqlb

\blemma\label{th3.2} Let $U_1= \inf\{t\ge 0: X_t\le Z_0(t)\}$. Under the assumptions of Theorem~\ref{th3.1}, we have
 \beqlb\label{eq3.13}
U_1= \inf\{t\ge 0: X_t= Z_0(t)\}= \inf\{t\ge 0: X_{t-}= Z_0(t-)\}.
 \eeqlb
\elemma

\begin{proof} From \eqref{eq3.8} and \eqref{eq3.12} it is easy to see that $X_{t-}= Z_0(t-)$ implies $X_t= Z_0(t)$. It follows that
 \beqnn
U_1\le \inf\{t\ge 0: X_t= Z_0(t)\}\le \inf\{t\ge 0: X_{t-}= Z_0(t-)\}.
 \eeqnn
Then \eqref{eq3.13} holds if $U_1= \infty$. In the case of $U_1< \infty$, we clearly have $X_{U_1-}\ge Z_0(U_1-)$ and $X_t> Z_0(t)$ for $0\le t< U_1$. By \eqref{eq3.8} and \eqref{eq3.12},
 \beqnn
\itDelta X_{U_1}\ar=\ar \int_{\{U_1\}}\int_0^\infty \int_0^{X_{U_1-}}\int_0^1 z M_0(\d s,\d z,\d u,\d v) \cr
 \ar\ar
+ \int_{\{U_1\}}\int_0^\infty \int_0^1 z N_0(\d s,\d z,\d v) \cr
 \ar\ge\ar
\int_{\{U_1\}}\int_0^\infty \int_0^{Z_0(U_1-)}\int_0^1 z M_0(\d s,\d z,\d u,\d v) \cr
 \ar\ar
+ \int_{\{U_1\}}\int_0^\infty \int_0^1 z N_0(\d s,\d z,\d v)= \itDelta Z_0(U_1).
 \eeqnn
On the other hand, by the right continuity of the processes,
 \beqnn\rv
X_{U_1-} + \itDelta X_{U_1}= X_{U_1}\le Z_0(U_1)= Z_0(U_1-) + \itDelta Z_0(U_1).
 \eeqnn
Then we must have $X_{U_1-}= Z_0(U_1-)$ and $\itDelta X_{U_1}= \itDelta Z_0(U_1)$, implying $X_{U_1}= Z_0(U_1)$. Those yield \eqref{eq3.13}. \end{proof}

By \eqref{eq3.13} and the definition of $U_1$ we have $X_{s-}> Z_0(s-)\ge 0$ for $0< s< U_1$. Let us consider the pure jump process $\{\xi_0(t): t\ge 0\}$ given by
 \beqlb\label{eq3.14}
~\qquad \xi_0(t)\ar=\ar \int_0^{t\land U_1}\int_0^\infty\int_0^{Z_0(s-)} \int_0^{\rho_1(Z_0(s-)-X_{s-},z)} [X_{s-}-Z_0(s-)] M_0(\d s,\d z,\d u,\d v) \\
 \ar\ar
+ \int_0^{t\land U_1}\int_0^\infty\int_0^{r_1(Z_0(s-)-X_{s-},z)} [X_{s-}-Z_0(s-)] N_0(\d s,\d z,\d v) \nnm\\
 \ar\ar
+ \int_0^{t\land U_1}\int_0^\infty\int_0^{Z_0(s-)} \int_{\rho_1(Z_0(s-)-X_{s-},z)} ^{\rho_2(X_{s-},Z_0(s-),z)} [Z_0(s-)-X_{s-}] M_0(\d s,\d z,\d u,\d v) \nnm\\
\ar\ar
+ \int_0^{t\land U_1}\int_0^\infty \int_{r_1(Z_0(s-)-X_{s-},z)}^{r_2(X_{s-},Z_0(s-),z)} [Z_0(s-)-X_{s-}] N_0(\d s,\d z,\d v). \nnm
 \eeqlb
The process makes its first jump at time $\sigma_1\land \tau_1$, where
 \beqnn
\sigma_1= \inf\{t\ge 0: \itDelta\xi_0(t)= X_{t-}-Z_0(t-)\}
 \eeqnn
and
 \beqnn
\tau_1= \inf\{t\ge 0: \itDelta\xi_0(t)= Z_0(t-)-X_{t-}\}.
 \eeqnn
In view of \eqref{eq3.3} and \eqref{eq3.4}, for any $x>y\ge 0$,
 \beqnn\rv
\int_0^\infty \rho_1(x-y,z)\mu(\d z)
 =
\int_0^\infty \rho_1(y-x,z)\mu(\d z)\le \frac{1}{2}\mu(x-y,\infty)< \infty
 \eeqnn
and
 \beqnn\rv
\int_0^\infty r_1(x-y,z)\nu(\d z)
 =
\int_0^\infty r_1(y-x,z)\nu(\d z)\le \frac{1}{2}\nu(x-y,\infty)< \infty.
 \eeqnn
Then the stopping times $\sigma_1$ and $\tau_1$ occur at the same rate
 \beqnn\rv
\bigg[Z_0(s-)\int_0^\infty \rho_1(X_{s-}-Z_0(s-),z)\mu(\d z) + \int_0^\infty r_1(X_{s-}-Z_0(s-),z)\nu(\d z)\bigg]\d s.
 \eeqnn
It follows that
 \beqlb\label{eq3.15}
\mbf{P}\big(\sigma_1< \tau_1\big|\sigma_1\land\tau_1< \infty\big)
 =
\mbf{P}\big(\tau_1< \sigma_1\big|\sigma_1\land\tau_1< \infty\big)= \frac{1}{2}.
 \eeqlb
For $t\ge 0$ let
 \beqlb\label{eq3.16}
Y_0(t)= Z_0(t\land U_1\land \sigma_1\land \tau_1) + \xi_0(t\land U_1\land \sigma_1\land \tau_1).
 \eeqlb

\blemma\label{th3.3} Under the assumptions of Theorem~\ref{th3.1}, we have $X_{\sigma_1}= Y_0(\sigma_1)> 0$ on the event $\{\sigma_1< \tau_1\le \infty\}$ and $X_{\tau_1}> Y_0(\tau_1)> 0$ on the event $\{\tau_1< \sigma_1\le \infty\}$. \elemma

\begin{proof} (1)~On the event $\{\sigma_1< \tau_1\le \infty\}$ we have $\sigma_1< U_1$ and so $\sigma_1< U_1\land \tau_1$. By \eqref{eq3.14} and the definition of $\sigma_1$ we have
 \beqlb\label{eq3.17}
\itDelta\xi_0(\sigma_1)= X_{\sigma_1-} - Z_0(\sigma_1-)> 0
 \eeqlb
and
 \beqnn
\ar\ar\int_{\{\sigma_1\}}\int_0^\infty \int_0^{Z_0(\sigma_1-)}\int_0^{\rho_1(Z_0(\sigma_1-)-X_{\sigma_1-},z)} M_0(\d s,\d z,\d u,\d v) \cr
 \ar\ar\qqquad
+ \int_{\{\sigma_1\}}\int_0^\infty \int_0^{r_1(Z_0(\sigma_1-)-X_{\sigma_1-},z)} N_0(\d s,\d z,\d v)= 1.
 \eeqnn
By the temporarily homogeneous nature of the Poisson random measures, we can enlarge the integration intervals of $u$ and $v$ without changing the equality above. In particular, we have
 \beqnn
\ar\ar\int_{\{\sigma_1\}}\int_0^\infty\int_0^{Z_0(\sigma_1-)}\int_0^1 M_0(\d s,\d z,\d u,\d v) + \int_{\{\sigma_1\}}\int_0^\infty \int_0^1 N_0(\d s,\d z,\d u) \cr
 \ar\ar\qquad
= \int_{\{\sigma_1\}}\int_0^\infty\int_0^{X_{\sigma_1-}}\int_0^1 M_0(\d s,\d z,\d u,\d v) + \int_{\{\sigma_1\}}\int_0^\infty \int_0^1 N_0(\d s,\d z,\d u)= 1.
 \eeqnn
It follows that
 \beqnn
z(\sigma_1)\ar:=\ar \int_{\{\sigma_1\}}\int_0^\infty\int_0^{Z_0(\sigma_1-)}\int_0^1 z M_0(\d s,\d z,\d u,\d v) + \int_{\{\sigma_1\}}\int_0^\infty \int_0^1 zN_0(\d s,\d z,\d u) \cr
 \ar=\ar
\int_{\{\sigma_1\}}\int_0^\infty\int_0^{X_{\sigma_1-}}\int_0^1 zM_0(\d s,\d z,\d u,\d v) + \int_{\{\sigma_1\}}\int_0^\infty \int_0^1 zN_0(\d s,\d z,\d u),
 \eeqnn
which means $z(\sigma_1)= \itDelta Z_0(\sigma_1)= \itDelta X_{\sigma_1}$ by \eqref{eq3.8} and \eqref{eq3.12}. This together with \eqref{eq3.17} yields
 \beqnn
Z_0(\sigma_1) + \xi_0(\sigma_1)= Z_0(\sigma_1-) + \itDelta Z_0(\sigma_1) + \itDelta\xi_0(\sigma_1)
 =
X_{\sigma_1-} + \itDelta X_{\sigma_1}= X_{\sigma_1}> 0.
 \eeqnn
By \eqref{eq3.16} we have $Y_0(\sigma_1)= X_{\sigma_1}> 0$.

(2)~On the event $\{\tau_1< \sigma_1\le \infty\}$ we have $\tau_1< U_1$ and so $\tau_1< U_1\land \sigma_1$. By \eqref{eq3.14} and the definition of $\tau_1$ we have
 \beqlb\label{eq3.18}
\itDelta\xi_0(\tau_1)= Z_0(\tau_1-)-X_{\tau_1-}< 0
 \eeqlb
and
 \beqnn
\ar\ar\int_{\{\tau_1\}}\int_0^\infty\int_0^{Z_0(\tau_1-)}\int_{\rho_1(Z_0(\tau_1-)-X_{\tau_1-},z)} ^{\rho_2(X_{\tau_1-},Z_0(\tau_1-),z)} M_0(\d s,\d z,\d u,\d v) \cr
 \ar\ar\qqquad\qquad
+ \int_{\{\tau_1\}}\int_0^\infty \int_{r_1(Z_0(\tau_1-)-X_{\tau_1-},z)}^{r_2(X_{\tau_1-},Z_0(\tau_1-),z)} N_0(\d s,\d z,\d v)= 1.
 \eeqnn
In view of \eqref{eq3.5} and \eqref{eq3.18}, for $0<z\le -\itDelta\xi_0(\tau_1)$ we have
 \beqnn
\rho_2(X_{\tau_1-},Z_0(\tau_1-),z)-\rho_1(Z_0(\tau_1-)-X_{\tau_1-},z)= \rho_1(-\itDelta\xi_0(\tau_1),z) =0
 \eeqnn
and
 \beqnn
r_2(X_{\tau_1-},Z_0(\tau_1-),z)-r_1(Z_0(\tau_1-)-X_{\tau_1-},z)= r_1(-\itDelta\xi_0(\tau_1),z) =0.
 \eeqnn
It follows that
 \beqlb\label{eq3.19}
z(\tau_1)\ar:=\ar \int_{\{\tau_1\}}\int_0^\infty\int_0^{Z_0(\tau_1-)}\int_{\rho_1(Z_0(\tau_1-)-X_{\tau_1-},z)} ^{\rho_2(X_{\tau_1-},Z_0(\tau_1-),z)} z M_0(\d s,\d z,\d u,\d v) \\
 \ar\ar
+ \int_{\{\tau_1\}}\int_0^\infty \int_{r_1(Z_0(\tau_1-)-X_{\tau_1-},z)} ^{r_2(X_{\tau_1-},Z_0(\tau_1-),z)} z N_0(\d s,\d z,\d u)> -\itDelta\xi_0(\tau_1)> 0. \nnm
 \eeqlb
As in Part~(1) of this the proof, one sees $z(\tau_1)= \itDelta Z_0(\tau_1)= \itDelta X_{\tau_1}$, and so
 \beqnn
0\le Z_0(\tau_1-)< Z_0(\tau_1-) + \itDelta Z_0(\tau_1) + \itDelta\xi_0(\tau_1)
 <
X_{\tau_1-} + \itDelta X_{\tau_1}= X_{\tau_1}.
 \eeqnn
Then we have $0< Y_0(\tau_1)< X_{\tau_1}$ by \eqref{eq3.16}. \end{proof}

\begin{proof}[Proof of Theorem~\ref{th3.1}] We give an inductive construction of the process $\{(Y_t,\xi_t): t\ge 0\}$ in a sequence of steps. The construction also yields the pathwise uniqueness of the solution of the equation system \eqref{eq3.9}-\eqref{eq3.10}.

\textit{Step~1.~} We start with the process $\{(X_t,Z_0(t)): t\ge 0\}$ defined by \eqref{eq3.8} and \eqref{eq3.12}. Consider separately the cases $\sigma_1= \tau_1= \infty$, $\sigma_1< \tau_1$ and $\tau_1< \sigma_1$. In the case of $\sigma_1= \tau_1= \infty$, we have $\xi_0(t)= 0$ for all $t\ge 0$. Then the process $\{(Y_t,\xi_t): t\ge 0\}$ is defined by $T=U_1$ and
 \beqnn
\xi_t= 0,~ Y_t= Z_0(t\land T) + X_t - X_{t\land T}.
 \eeqnn
In the case of $\sigma_1< \tau_1\le \infty$, we have $X_{\sigma_1}= Y_0(\sigma_1)> 0$ by Lemma~\ref{th3.3}. Then the process $\{(Y_t,\xi_t): t\ge 0\}$ is given by $T= \sigma_1$ and
 \beqnn
\xi_t= \xi_0(t\land T), ~ Y_t= Y_0(t\land T) + X_t - X_{t\land T}.
 \eeqnn
In case of $\tau_1< \sigma_1\le \infty$, we clearly have $Y_t= Z_0(t)$ and $\xi_t= 0$ for $0\le t< \tau_1$. The continuing construction of $\{(Y_{\tau_1+t},\xi_{\tau_1+t}): t\ge 0\}$ is given in the next step.

\textit{Step~2.~} Suppose that $\tau_1< \sigma_1\le \infty$. In this case, we have $X_{\tau_1}> Y_0(\tau_1)> 0$ by Lemma~\ref{th3.3}. Let $X_1(t)= X_{\tau_1+t}$ for $t\ge 0$. Then $\{X_1(t): t\ge 0\}$ is also a CBIC-process. From \eqref{eq3.8} it follows that
 \beqlb\label{eq3.20}
X_1(t)= X_{\tau_1} + \int_0^t\int_0^{X_1(s-)}L_1(\d s,\d u) - \int_0^t g(X_1(s))\d s + \eta_{\tau_1+t},
 \eeqlb
where $L_1(\d s,\d u)= L_0(\tau_1+\d s,\d u)$. By Theorem~\ref{th2.1} we can construct another CBIC-process $\{Z_1(t): t\ge 0\}$ by the pathwise unique solution to
 \beqlb\label{eq3.21}
Z_1(t)= Y_0(\tau_1) + \int_0^t\int_0^{Z_1(s-)} L_1^*(\d s,\d u) - \int_0^t g(Z_1(s)) \d s + \eta_{\tau_1+t},
 \eeqlb
where $L_1^*(\d s,\d u)= L_0^*(\tau_1+\d s,\d u)$. Then we repeat the procedure of in the first step for the process $\{(X_1(t),Z_1(t)): t\ge 0\}$.

In view of \eqref{eq3.15}, we only need \textit{a finite number} of sequential steps to complete the construction of $\{(Y_t,\xi_t): t\ge 0\}$, which means that the process $\{\xi_t: t\ge 0\}$ has at most a finite number of jumps. By the construction, we have \eqref{eq3.11} and $X_{T+t}= Y_{T+t}$ for every $t\ge 0$. \end{proof}

\bremark\label{th3.4} By Theorem~\ref{th3.1}, the process $\{(X_t,Y_t): t\ge 0\}$ defined by \eqref{eq3.8} and \eqref{eq3.9} with $X_0\ge Y_0\ge 0$ actually lives in the space $D:= \{(x,y): x\ge y\ge 0\}\subset \mbb{R}_+^2$. Here we think of the coupling time $T$ as a part of the solution. We may also say that $\{(X_t,Y_t,\xi_t): t\ge 0\}$ is a pathwise unique solution to the system of equations \eqref{eq3.8}, \eqref{eq3.9} and \eqref{eq3.10}. \eremark

\btheorem\label{th3.5} The solution $\{Y_t: t\ge 0\}$ to \eqref{eq3.9} is a CBIC-process. \etheorem

\begin{proof} Let $\{N(\d s,\d z)\}$ be the optional random measure defined by, for $t\ge 0$ and $A\in \mcr{B}(0,\infty)$,
 \beqnn
N((0,t]\times A)\ar=\ar \int_0^t\int_0^\infty\int_0^{r_1(Y_{s-}-X_{s-},z)} \I_A(z+X_{s-}-Y_{s-}) N_0(\d s,\d z,\d v) \cr
 \ar\ar
+ \int_0^t\int_0^\infty\int_{r_1(Y_{s-}-X_{s-},z)}^{r_2(X_{s-},Y_{s-},z)} \I_A(z+Y_{s-}-X_{s-}) N_0(\d s,\d z,\d v) \cr
 \ar\ar
+ \int_0^t\int_0^\infty\int_{r_2(X_{s-},Y_{s-},z)}^1 \I_A(z) N_0(\d s,\d z,\d v).
 \eeqnn
Then $N(\d s,\d z)$ has predictable compensator $\hat{N}(\d s,\d z)$ determined by
 \beqnn
\hat{N}((0,t]\times A)\ar=\ar \int_0^t\d s\int_0^\infty \I_A(z+X_{s-}-Y_{s-})r_1(Y_{s-}-X_{s-},z) \nu(\d z) \cr
 \ar\ar
+ \int_0^t\d s\int_0^\infty \I_A(z+Y_{s-}-X_{s-})r_1(X_{s-}-Y_{s-},z) \nu(\d z) \cr
 \ar\ar
+ \int_0^t\d s\int_0^\infty \I_A(z)[1-r_2(X_{s-},Y_{s-},z)]\nu(\d z) \cr
 \ar=\ar
\int_0^t\d s\int_0^\infty \I_A(z+X_{s-}-Y_{s-})\nu(Y_{s-}-X_{s-},\d z) \cr
 \ar\ar
+ \int_0^t\d s\int_0^\infty \I_A(z+Y_{s-}-X_{s-})\nu(X_{s-}-Y_{s-},\d z) \cr
 \ar\ar
+ \int_0^t\d s\int_0^\infty \I_A(z)[\nu(\d z)-\nu(X_{s-},Y_{s-},\d z)].
 \eeqnn
In view of \eqref{eq3.2}, for any $x\ge y\ge 0$ we have
 \beqnn
\ar\ar\int_0^\infty \I_A(z+x-y) \nu(y-x,\d z) + \int_0^\infty \I_A(z+y-x) \nu(x-y,\d z) \cr
 \ar\ar\qquad
=\int_0^\infty \I_A(z) \nu(x-y,\d z) + \int_0^\infty \I_A(z) \nu(y-x,\d z) \cr
 \ar\ar\qquad
=\int_0^\infty \I_A(z) \nu(x,y,\d z).
 \eeqnn
Then $\hat{N}((0,t]\times A)= t\nu(A)$, and so $N(\d s,\d z)$ is a Poisson random measure on $(0,\infty)^2$ with intensity $\d s\nu(\d z)$; see, e.g., Theorem~III.6.2 in Ikeda and Watanabe \cite[p.75]{IkW89}. Let $\{M(\d s,\d z,\d u)\}$ be the optional random measure defined by, for $t\ge 0$ and $A,B\in \mcr{B}(0,\infty)$,
 \beqnn
M((0,t]\times A\times B)\ar=\ar \int_0^t\int_0^\infty\int_B\int_0^{\rho_1(Y_{s-}-X_{s-},z)} \I_A(z+X_{s-}-Y_{s-}) M_0(\d s,\d z,\d u,\d v) \cr
 \ar\ar
+ \int_0^t\int_0^\infty\int_B\int_{\rho_1(Y_{s-}-X_{s-},z)}^{\rho_2(X_{s-},Y_{s-},z)} \I_A(z+Y_{s-}-X_{s-}) M_0(\d s,\d z,\d u,\d v) \cr
 \ar\ar
+ \int_0^t\int_0^\infty\int_B\int_{\rho_2(X_{s-},Y_{s-},z)}^1 \I_A(z) M_0(\d s,\d z,\d u,\d v).
 \eeqnn
By similar {\rrv calculations as above} one can see $\{M(\d s,\d z,\d u)\}$ is a Poisson random measure on $(0,\infty)^3$ with intensity $\d s\mu(\d z)\d u$. Let $\{\eta(t)\}$ be the $(\mcr{F}_t)$-subordinator defined by \eqref{eq1.11}. Let $\{L(\d s,\d u)\}$ be the time-space $(\mcr{F}_t)$-L\'evy noise defined by
 \beqnn
L(\d s,\d u)\ar=\ar -\, \I_{\{s\le T\}}W_0(\d s,\d u) - b\d s\d u + \int_{\{0< z\le 1\}} z \tilde{M}(\d s,\d z,\d u) \cr
 \ar\ar
+\, \I_{\{s>T\}}W_0(\d s,\d u) + \int_{\{1< z< \infty\}} z M(\d s,\d z,\d u)
 \eeqnn
From \eqref{eq3.9} and \eqref{eq3.10} we obtain
 \beqnn
Y_t= Y_0 + \int_0^t\int_0^{Y_{s-}} L(\d s,\d u) - \int_0^t g(Y_s) \d s + \eta(t).
 \eeqnn
Then $\{Y_t: t\ge 0\}$ is a CBIC-process by the uniqueness of the solution to the equation. \end{proof}

The pathwise uniqueness of the solutions to \eqref{eq3.8} and \eqref{eq3.9} implies that $\{(X_t,Y_t): t\ge 0\}$ is a Markov process. Then it is a Markov coupling of the CBIC-process. Since zero is a trap for $\{X_t-Y_t: t\ge 0\}$, we have
 \beqlb\label{eq3.22}
\quad X_t-Y_t\ar=\ar X_0-Y_0 + \int_0^t\int_{Y_s}^{X_s} W_0(\d s,\d u) + 2\int_0^{t\wedge T}\int_0^{Y_s} W_0(\d s,\d u) \\
 \ar\ar
- \int_0^t \big[g(X_s)-g(Y_s)\big]\d s + \int_0^t\int_0^1\int_{Y_{s-}}^{X_{s-}}\int_0^1 z \tilde{M}_0(\d s,\d z,\d u,\d v) \nnm\\
 \ar\ar
+ \int_0^t\int_1^\infty\int_{Y_{s-}}^{X_{s-}}\int_0^1 z M_0(\d s,\d z,\d u,\d v) - \xi_t. \nnm
 \eeqlb

\bremark\label{th3.6} Concerning the coupling process, what is important to us is its movement before the succeeding time. To establish \eqref{eq1.19} and \eqref{eq1.31}, we need to make the coupling succeed as early as possible. The Gaussian and Poissonian integrals on the right-hand side of \eqref{eq3.22} are important in bringing the process to the success. Clearly, the second Gaussian integral comes from the reflection of the Gaussian component in \eqref{eq3.7}. The Poisson noises cannot be reflected directly as the CBIC-processes have no negative jumps. The pure jump process $\{\xi_t: t\ge 0\}$ defined by \eqref{eq3.10} gives a proper formulation of the reflections of the discontinuous noises. In view of \eqref{eq3.15}, a possible jump of this process at time $s>0$ would take the values $X_{s-}-Y_{s-}> 0$ and $Y_{s-}-X_{s-}< 0$ with equal probability $1/2$. The introduction of the process $\{\xi_t: t\ge 0\}$ was inspired by the construction of a coupling process of Luo and Wang \cite{LuW19}; see also \cite{LiW20, LMW21}. In fact, if Grey's condition \eqref{eq1.9} is not satisfied, the success of the coupling process can only come with the first positive jump of $\{\xi_t: t\ge 0\}$. \eremark

\section{The coupling generator}

 \setcounter{equation}{0}

{\rv In this section, we assume the CBIC-process with an arbitrary initial value is conservative. We shall give a characterization for the generator of the Markov coupling process $\{(X_t,Y_t): t\ge 0\}$ defined by \eqref{eq3.8} and \eqref{eq3.9} in terms of a martingale problem.} We also establish some estimates for the generator, which provide the basis for the proof of the exponential ergodicity. Recall that $\{(X_t,Y_t): t\ge 0\}$ has state space $D:= \{(x,y): x\ge y\ge 0\}$. Let $\itDelta= \{(z,z): z\ge 0\}\subset D$ and $\itDelta^c= D\setminus \itDelta$. Given a function $F$ on $D$ twice continuously differentiable on $\itDelta^c$, we write
 \beqlb\label{eq4.1}
\tilde{L}F(x,y)= \tilde{L}_0F(x,y) + \tilde{L}_1F(x,y), \quad (x,y)\in \itDelta^c,
 \eeqlb
where
 \beqlb\label{eq4.2}
~\qquad \tilde{L}_0F(x,y)\ar=\ar \big[\beta-bx-g(x)\big]F'_x(x,y) + \big[\beta-by-g(y)\big]F'_y(x,y) + cxF''_{xx}(x,y) \\
 \ar\ar
+\,cyF''_{yy}(x,y) - 2cyF''_{xy}(x,y) + \int_0^\infty \big[F(x+z,y+z)-F(x,y)\big] \nu(\d z) \nnm\\
 \ar\ar
+\, (x-y)\int_0^\infty\big[F(x+z,y)-F(x,y)-F_x'(x,y)z\I_{\{z\le 1\}}\big] \mu(\d z) \nnm\\
 \ar\ar
+\,y\int_0^\infty \big[F(x+z,y+z)-F(x,y)-(F_x'+F_y')(x,y)z\I_{\{z\le 1\}}\big] \mu(\d z) \nnm
 \eeqlb
and
 \beqlb\label{eq4.3}
\tilde{L}_1F(x,y)\ar=\ar y\int_0^\infty \big[F(x+z,2y+z-x)-F(x+z,y+z)\big] \mu_{x-y}(\d z) \\
 \ar\ar
+\,y\int_0^\infty \big[F(x+z,x+z)-F(x+z,y+z)\big] \mu_{y-x}(\d z) \nnm\\
 \ar\ar
+ \int_0^\infty \big[F(x+z,2y+z-x)-F(x+z,y+z)\big] \nu_{x-y}(\d z) \nnm\\
 \ar\ar
+ \int_0^\infty \big[F(x+z,x+z)-F(x+z,y+z)\big] \nu_{y-x}(\d z). \nnm
 \eeqlb
Let $\mcr{D}(\tilde{L})$ denote the linear space consisting of the functions $F$ such that the integrals in \eqref{eq4.2} and \eqref{eq4.3} are convergent and define functions locally bounded on compact subsets of $\itDelta^c$. The operator $\tilde{L}$ determines the movement of the coupling process before its succeeding time. In particular, the component $\tilde{L}_1$ is induced by the process $\{\xi_t: t\ge 0\}$ defined by \eqref{eq3.10}. We call $(\tilde{L},\mcr{D}(\tilde{L}))$ the \textit{coupling generator} of the CBIC-process. The precise meaning of this terminology is made clear by the martingale problem given in the following:

\btheorem\label{th4.1} Let $\{(X_t,Y_t): t\ge 0\}$ be the Markov coupling process defined by \eqref{eq3.8} and \eqref{eq3.9} with $X_0> Y_0\ge 0$ and let $\zeta_n^*= \inf\{t\ge 0: X_t\ge n$ or $X_t-Y_t\le 1/n\}$. Then for any $n\ge 1$ and $F\in \mcr{D}(\tilde{L})$ we have
 \beqlb\label{eq4.4}
F(X_{t\land \zeta_n^*},Y_{t\land \zeta_n^*}) = F(X_0,Y_0) + \int_0^{t\land \zeta_n^*} \tilde{L}F(X_s,Y_s) \d s + M_n(t),
 \eeqlb
where $\{M_n(t): t\ge 0\}$ is a martingale.
\etheorem

\begin{proof} By \eqref{eq3.8}, \eqref{eq3.9} and \eqref{eq3.10}, the process $\{(X_{t\land \zeta_n^*},Y_{t\land \zeta_n^*}): t\ge 0\}$ is a semimartingale taking values in $\itDelta\cup D_n$, where $D_n= \{(x,y)\in D: x-y\ge 1/n\}$. In view of the three equations, a possible jump $(\itDelta X_s,\itDelta Y_s)$ of the process at time $s> 0$ is brought about by a point $(s,z,u,v)\in \supp(M_0)$ with $0< u\le X_{s-}$ or by a point $(s,z,v)\in \supp(N_0)$. The details of the jump in the two cases are given, respectively, by
 \beqnn
\quad(\itDelta X_s,\itDelta Y_s)= \left\{\begin{array}{ll}
 (z,z+X_{s-}-Y_{s-}), & u\in (0,Y_{s-}], v\in (0,\hat{\rho}_1(s,z)], \smallskip\cr
 (z,z+Y_{s-}-X_{s-}), & u\in (0,Y_{s-}], v\in (\hat{\rho}_1(s,z),\hat{\rho}_2(s,z)], \smallskip\cr
 (z,z), & u\in (0,Y_{s-}], v\in (\hat{\rho}_2(s,z),1], \smallskip\cr
 (z,0), & u\in (Y_{s-},X_{s-}], v\in (0,1],
\end{array}\right.
 \eeqnn
where $\hat{\rho}_1(s,z)= \rho_1(X_{s-}-Y_{s-},z)$ and $\hat{\rho}_2(s,z)= \rho_2(X_{s-},Y_{s-},z)$, and
 \beqnn
~~(\itDelta X_s,\itDelta Y_s)= \left\{\begin{array}{ll}
 (z,z+X_{s-}-Y_{s-}), & v\in (0,\hat{r}_1(s,z)], \smallskip\cr
 (z,z+Y_{s-}-X_{s-}), & v\in (\hat{r}_1(s,z),\hat{r}_2(s,z)], \smallskip\cr
 (z,z), & v\in (\hat{r}_2(s,z),1].
\end{array}\right.
 \eeqnn
where $\hat{r}_1(s,z)= r_1(X_{s-}-Y_{s-},z)$ and $\hat{r}_2(s,z)= r_2(X_{s-},Y_{s-},z)$. Then the coupling process may have totally seven different types of jumps. Observe also that $(X_{s-},Y_{s-})\in D_n\subset \itDelta^c$ when $0< s\le \zeta_n^*$. For any $(x,y)\in \itDelta^c$ write
 \beqnn
\tilde{L}_2F(x,y) = cxF''_{xx}(x,y) + cyF''_{yy}(x,y) - 2cyF''_{xy}(x,y) - g(x)F'_x(x,y) - g(y)F'_y(x,y).
 \eeqnn
We can use It\^o's formula to see
 \beqnn
F(X_{t\land \zeta_n^*},Y_{t\land \zeta_n^*})\ar=\ar F(X_0,Y_0) + \int_0^{t\land \zeta_n^*} \tilde{L}_2F(X_s,Y_s)\d s + {\rv \int_0^{t\land \zeta_n^*}\int_0^{X_{s-}} F'_x(X_{s-},Y_{s-}) L_0(\d s,\d u)} \cr
 \ar\ar
+ {\rv \int_0^{t\land \zeta_n^*}\int_0^{Y_{s-}} F'_y(X_{s-},Y_{s-}) L^*_0(\d s,\d u)} + \int_0^{t\land \zeta_n^*} (F'_x+F'_y)(X_{s-},Y_{s-}) \d\eta_s \cr
 \ar\ar
+\int_0^{t\land \zeta_n^*}\int_0^\infty\int_0^{Y_{s-}}\int_0^{\hat{\rho}_1(s,z)} F'_y(X_{s-},Y_{s-}) (X_{s-}-Y_{s-}) M_0(\d s,\d z,\d u,\d v) \cr
 \ar\ar
+ \int_0^{t\land \zeta_n^*}\int_0^\infty\int_0^{\hat{r}_1(s,z)} F'_y(X_{s-},Y_{s-}) (X_{s-}-Y_{s-}) N_0(\d s,\d z,\d v) \cr
 \ar\ar
+ \int_0^{t\land\zeta_n^*}\int_0^\infty\int_0^{Y_{s-}}\int_{\hat{\rho}_1(s,z)}^{\hat{\rho}_2(s,z)} F'_y(X_{s-},Y_{s-})(Y_{s-}-X_{s-}) M_0(\d s,\d z,\d u,\d v) \cr
 \ar\ar
+ \int_0^{t\land \zeta_n^*}\int_0^\infty\int_{\hat{r}_1(s,z)}^{\hat{r}_2(s,z)} F'_y(X_{s-},Y_{s-}) (Y_{s-}-X_{s-}) N_0(\d s,\d z,\d v) \cr
 \ar\ar
+ \int_0^{t\land \zeta_n^*}\int_0^\infty\int_0^{Y_{s-}}\int_0^{\hat{\rho}_1(s,z)} \big[F(X_{s-}+z,X_{s-}+z) - F(X_{s-},Y_{s-}) \nnm\\
 \ar\ar\quad
-\, zF'_x(X_{s-},Y_{s-}) - (z+X_{s-}-Y_{s-})F'_y(X_{s-},Y_{s-})\big] M_0(\d s,\d z,\d u,\d v) \nnm\\
 \ar\ar
+ \int_0^{t\land \zeta_n^*}\int_0^\infty\int_0^{Y_{s-}}\int_{\hat{\rho}_1(s,z)}^{\hat{\rho}_2(s,z)} \big[F(X_{s-}+z,2Y_{s-}+z-X_{s-}) - F(X_{s-},Y_{s-}) \nnm\\
 \ar\ar\quad
-\, zF'_x(X_{s-},Y_{s-}) - (z+Y_{s-}-X_{s-})F'_y(X_{s-},Y_{s-})\big] M_0(\d s,\d z,\d u,\d v) \nnm\\
 \ar\ar
+ \int_0^{t\land \zeta_n^*}\int_0^\infty\int_0^{Y_{s-}}\int_{\hat{\rho}_2(s,z)}^1 \big[F(X_{s-}+z,Y_{s-}+z) - F(X_{s-},Y_{s-}) \nnm\\
 \ar\ar\qqquad\qqquad
-\,z(F'_x+F'_y)(X_{s-},Y_{s-})\big] M_0(\d s,\d z,\d u,\d v) \nnm\\
 \ar\ar
+ \int_0^{t\land \zeta_n^*}\int_0^\infty\int_{Y_{s-}}^{X_{s-}} \int_0^1 \big[F(X_{s-}+z,Y_{s-}) - F(X_{s-},Y_{s-}) \nnm\\
 \ar\ar\qqquad\qqquad
-\,zF'_x(X_{s-},Y_{s-})\big] M_0(\d s,\d z,\d u,\d v) \cr
 \ar\ar
+ \int_0^{t\land \zeta_n^*}\int_0^\infty\int_0^{\hat{r}_1(s,z)} \big[F(X_{s-}+z,X_{s-}+z) - F(X_{s-},Y_{s-}) \nnm\\
 \ar\ar\qquad
-\,zF'_x(X_{s-},Y_{s-}) - (z+X_{s-}-Y_{s-})F'_y(X_{s-},Y_{s-})\big] N_0(\d s,\d z,\d v) \nnm\\
 \ar\ar
+ \int_0^{t\land \zeta_n^*}\int_0^\infty\int_{\hat{r}_1(s,z)}^{\hat{r}_2(s,z)} \big[F(X_{s-}+z,2Y_{s-}+z-X_{s-}) - F(X_{s-},Y_{s-}) \nnm\\
 \ar\ar\qquad
-\, zF'_x(X_{s-},Y_{s-}) - (z+Y_{s-}-X_{s-})F'_y(X_{s-},Y_{s-})\big] N_0(\d s,\d z,\d v) \nnm\\
 \ar\ar
+ \int_0^{t\land \zeta_n^*}\int_0^\infty\int_{\hat{r}_2(s,z)}^1 \big[F(X_{s-}+z,Y_{s-}+z) - F(X_{s-},Y_{s-}) \nnm\\
 \ar\ar\qqquad\qqquad
-\,z(F'_x+F'_y)(X_{s-},Y_{s-})\big] N_0(\d s,\d z,\d v) \nnm\\
 \ar=\ar
F(X_0,Y_0) + \int_0^{t\land \zeta_n^*} \tilde{L}_2F(X_s,Y_s)\d s + {\rv \int_0^{t\land \zeta_n^*}\int_0^{X_{s-}} F'_x(X_{s-},Y_{s-}) L_0(\d s,\d u)} \cr
 \ar\ar
+ {\rv \int_0^{t\land \zeta_n^*}\int_0^{Y_{s-}} F'_y(X_{s-},Y_{s-}) L^*_0(\d s,\d u)} + \int_0^{t\land \zeta_n^*} (F'_x+F'_y)(X_{s-},Y_{s-}) \d\eta_s \cr
 \ar\ar
+ \int_0^{t\land \zeta_n^*}\int_0^\infty\int_0^{Y_{s-}}\int_0^{\hat{\rho}_1(s,z)} \big[F(X_{s-}+z,X_{s-}+z) - F(X_{s-},Y_{s-}) \nnm\\
 \ar\ar\qqquad\qqquad
-\,z(F'_x+F'_y)(X_{s-},Y_{s-})\big] M_0(\d s,\d z,\d u,\d v) \nnm\\
 \ar\ar
+ \int_0^{t\land \zeta_n^*}\int_0^\infty\int_0^{Y_{s-}} \int_{\hat{\rho}_1(s,z)}^{\hat{\rho}_2(s,z)} \big[F(X_{s-}+z,2Y_{s-}+z-X_{s-}) \nnm\\
 \ar\ar\qqquad
-\, F(X_{s-},Y_{s-}) - z(F'_x+F'_y)(X_{s-},Y_{s-})\big] M_0(\d s,\d z,\d u,\d v) \nnm\\
 \ar\ar
+ \int_0^{t\land \zeta_n^*}\int_0^\infty\int_0^{Y_{s-}}\int_{\hat{\rho}_2(s,z)}^1 \big[F(X_{s-}+z,Y_{s-}+z) - F(X_{s-},Y_{s-}) \nnm\\
 \ar\ar\qqquad\qqquad
-\,z(F'_x+F'_y)(X_{s-},Y_{s-})\big] M_0(\d s,\d z,\d u,\d v) \nnm\\
 \ar\ar
+ \int_0^{t\land \zeta_n^*}\int_0^\infty\int_{Y_{s-}}^{X_{s-}} \int_0^1 \big[F(X_{s-}+z,Y_{s-}) - F(X_{s-},Y_{s-}) \nnm\\
 \ar\ar\qqquad\qqquad
-\,zF'_x(X_{s-},Y_{s-})\big] M_0(\d s,\d z,\d u,\d v) \nnm\\
 \ar\ar
+ \int_0^{t\land \zeta_n^*}\int_0^\infty\int_0^{\hat{r}_1(s,z)} \big[F(X_{s-}+z,X_{s-}+z) - F(X_{s-},Y_{s-}) \nnm\\
 \ar\ar\qqquad\qqquad
-\,z(F'_x+F'_y)(X_{s-},Y_{s-})\big] N_0(\d s,\d z,\d v) \nnm\\
 \ar\ar
+ \int_0^{t\land \zeta_n^*}\int_0^\infty\int_{\hat{r}_1(s,z)}^{\hat{r}_2(s,z)} \big[F(X_{s-}+z,2Y_{s-}+z-X_{s-}) \nnm\\
 \ar\ar\qqquad
-\, F(X_{s-},Y_{s-}) - z(F'_x+F'_y)(X_{s-},Y_{s-})\big] N_0(\d s,\d z,\d v) \nnm\\
 \ar\ar
+ \int_0^{t\land \zeta_n^*}\int_0^\infty\int_{\hat{r}_2(s,z)}^1 \big[F(X_{s-}+z,Y_{s-}+z) - F(X_{s-},Y_{s-}) \nnm\\
 \ar\ar\qqquad\qqquad
-\,z(F'_x+F'_y)(X_{s-},Y_{s-})\big] N_0(\d s,\d z,\d v).
 \eeqnn
By some further cancellations of the terms,
 \beqlb\label{eq4.5}
~\quad F(X_{t\land \zeta_n^*},Y_{t\land \zeta_n^*})
 \ar=\ar
F(X_0,Y_0) + \int_0^{t\land \zeta_n^*} \tilde{L}_2F(X_{s},Y_{s})\d s + \int_0^{t\land \zeta_n^*} (\beta-bX_s)F'_x(X_s,Y_s) \d s \\
 \ar\ar
+ \int_0^{t\land \zeta_n^*}\int_0^{X_s} F'_x(X_s,Y_s) W_0(\d s,\d u) \nnm\\
 \ar\ar
+ \int_0^{t\land \zeta_n^*}\int_0^{X_{s-}}\int_0^1\int_0^1 F'_x(X_{s-},Y_{s-})z \tilde{M}_0(\d s,\d z,\d u,\d v) \nnm\\
 \ar\ar
+ \int_0^{t\land \zeta_n^*} (\beta-bY_s)F'_y(X_s,Y_s) \d s - \int_0^{t\land \zeta_n^*}\int_0^{Y_s} F'_y(X_s,Y_s) W_0(\d s,\d u) \nnm\\
\ar\ar
+ \int_0^{t\land \zeta_n^*}\int_0^{Y_{s-}}\int_0^1\int_0^1 F'_y(X_{s-},Y_{s-})z \tilde{M}_0(\d s,\d z,\d u,\d v) \nnm\\
 \ar\ar
+ \int_0^{t\land \zeta_n^*}\int_0^\infty\int_0^{Y_{s-}}\int_0^{\hat{\rho}_1(s,z)} \big[F(X_{s-}+z,X_{s-}+z) - F(X_{s-},Y_{s-}) \nnm\\
 \ar\ar\qqquad\qquad
-\,\I_{\{z\le 1\}}z(F'_x+F'_y)(X_{s-},Y_{s-})\big] \tilde{M}_0(\d s,\d z,\d u,\d v) \nnm\\
 \ar\ar
+ \int_0^{t\land \zeta_n^*}\int_0^\infty\int_0^{Y_{s-}}\int_{\hat{\rho}_1(s,z)}^{\hat{\rho}_2(s,z)} \big[F(X_{s-}+z,2Y_{s-}+z-X_{s-}) \nnm\\
 \ar\ar\quad
-\, F(X_{s-},Y_{s-}) - \I_{\{z\le 1\}}z(F'_x+F'_y)(X_{s-},Y_{s-})\big] \tilde{M}_0(\d s,\d z,\d u,\d v) \nnm\\
 \ar\ar
+ \int_0^{t\land \zeta_n^*}\int_0^\infty\int_0^{Y_{s-}}\int_{\hat{\rho}_2(s,z)}^1 \big[F(X_{s-}+z,Y_{s-}+z) - F(X_{s-},Y_{s-}) \nnm\\
 \ar\ar\qqquad\qquad
-\,\I_{\{z\le 1\}}z(F'_x+F'_y)(X_{s-},Y_{s-})\big)\big] \tilde{M}_0(\d s,\d z,\d u,\d v) \nnm\\
 \ar\ar
+ \int_0^{t\land \zeta_n^*}\int_0^\infty\int_{Y_{s-}}^{X_{s-}} \int_0^1 \big[F(X_{s-}+z,Y_{s-}) - F(X_{s-},Y_{s-}) \nnm\\
 \ar\ar\qqquad\qqquad
-\,\I_{\{z\le 1\}}zF'_x(X_{s-},Y_{s-})\big] \tilde{M}_0(\d s,\d z,\d u,\d v) \nnm\\
 \ar\ar
+ \int_0^{t\land \zeta_n^*}\int_0^\infty\int_0^{\hat{r}_1(s,z)} \big[F(X_{s-}+z,X_{s-}+z) \nnm\\
 \ar\ar\qqquad\qqquad
-\, F(X_{s-},Y_{s-})\big] \tilde{N}_0(\d s,\d z,\d v) \nnm\\
 \ar\ar
+ \int_0^{t\land \zeta_n^*}\int_0^\infty\int_{\hat{r}_1(s,z)}^{\hat{r}_2(s,z)} \big[F(X_{s-}+z,2Y_{s-}+z-X_{s-}) \nnm\\
 \ar\ar\qqquad\qqquad
-\, F(X_{s-},Y_{s-})\big] \tilde{N}_0(\d s,\d z,\d v) \nnm\\
 \ar\ar
+ \int_0^{t\land \zeta_n^*}\int_0^\infty\int_{\hat{r}_2(s,z)}^1 \big[F(X_{s-}+z,Y_{s-}+z) \nnm\\
 \ar\ar\qqquad\qqquad
-\,F(X_{s-},Y_{s-})\big] \tilde{N}_0(\d s,\d z,\d v) \nnm\\
 \ar\ar
+ \int_0^{t\land \zeta_n} \tilde{L}_3F(X_{s},Y_{s})\d s, \nnm
 \eeqlb
where
 \beqnn
\tilde{L}_3F(x,y)
 \ar=\ar
y\int_0^\infty \big[F(x+z,y+z)-F(x,y)-(F_x'+F_y')(x,y)z\I_{\{z\le 1\}}\big] \mu(\d z) \cr
 \ar\ar
+\,y\int_0^\infty \big[F(x+z,x+z)-F(x+z,y+z)\big] \mu_{y-x}(\d z) \cr
 \ar\ar
+\, y\int_0^\infty \big[F(x+z,2y+z-x)-F(x+z,y+z)\big] \mu_{x-y}(\d z) \cr
 \ar\ar
+\,(x-y)\int_0^\infty\big[F(x+z,y)-F(x,y)-F_x'(x,y)z\I_{\{z\le 1\}}\big] \mu(\d z) \cr
 \ar\ar
+\int_0^\infty \big[F(x+z,y+z)-F(x,y)\big] \nu(\d z) \cr
 \ar\ar
+ \int_0^\infty \big[F(x+z,x+z)-F(x+z,y+z)\big] \nu_{y-x}(\d z) \cr
 \ar\ar
+ \int_0^\infty \big[F(x+z,2y+z-x)-F(x+z,y+z)\big] \nu_{x-y}(\d z).
 \eeqnn
Clearly, we have
 \beqnn
\tilde{L}F(x,y)= (\beta-bx)F'_x(x,y) + (\beta-by)F'_y(x,y) + \tilde{L}_2F(x,y) + \tilde{L}_3F(x,y).
 \eeqnn
Since $Y_{s-}\le X_{s-}\le n$ for $0<s\le \zeta_n^*$, the stochastic integrals on the right-hand side of \eqref{eq4.5} define a martingale $\{M_n(t): t\ge 0\}$. Then the desired result follows. \end{proof}

We next consider {\rv two special forms of the function} $F\in \mcr{D}(\tilde{L})$. Recall that $C_b^2(\mbb{R}_+)$ denotes the space of bounded and continuous functions on $\mbb{R}_+$ with bounded and continuous derivatives up to the second order.

{\rv

\blemma\label{th4.2'} Suppose that $f\in C_b^2(\mbb{R}_+)$ is a nonnegative, nondecreasing and concave function. Let
 \beqlb\label{eq4.6'}
F(x,y)= f(x-y)\I_{\{x\neq y\}}, \quad (x,y)\in D.
 \eeqlb
Then $F\in \mcr{D}(\tilde{L})$ and, for $(x,y)\in \itDelta^c$,
 \beqlb\label{eq4.7'}
\tilde{L}F(x,y)\ar\le\ar (x-y)\bigg[cf''(x-y) - bf'(x-y) + \int_0^\infty \big[f(x-y+z) - f(x-y) \\
 \ar\ar\qqquad\qqquad
-\,f'(x-y)z\I_{\{z\le 1\}}\big] \mu(\d z)\bigg]. \nnm
 \eeqlb
\elemma

\begin{proof} By \eqref{eq4.2} and \eqref{eq4.3} it is clear that $F\in \mcr{D}(\tilde{L})$. Since $g$ is nondecreasing and $f$ is nonnegative, nondecreasing and concave, by \eqref{eq4.2} we have, for $(x,y)\in \itDelta^c$,
 \beqnn
\tilde{L}_0F(x,y)\ar=\ar - \big[bx-by+g(x)-g(y)\big]f'(x-y) + c(x+3y)f''(x-y) \nnm\\
 \ar\ar
+\,(x-y)\int_0^\infty \big[f(x-y+z)-f(x-y) - f'(x-y)z\I_{\{z\le 1\}}\big] \mu(\d z) \nnm\\
 \ar\le\ar
(x-y)\bigg[cf''(x-y) - bf'(x-y) + \int_0^\infty \big[f(x-y+z) - f(x-y) \\
 \ar\ar\qqquad
-\,f'(x-y)z\I_{\{z\le 1\}}\big] \mu(\d z)\bigg].
\eeqnn
On the other hand, since $f(2z)-2f(z)\le0$ for any $z\ge 0$, by \eqref{eq4.3} it is easy to see that
\beqnn
\tilde{L}_1F(x,y)= y\big[f(2(x-y)) - 2f(x-y)\big]\mu_{x-y}(0,\infty) + \big[f(2(x-y)) - 2f(x-y)\big]\nu_{x-y}(0,\infty)\le 0.
 \eeqnn
These together with \eqref{eq4.1} yield \eqref{eq4.7'}. \end{proof}

}

\blemma\label{th4.2} Suppose that $f\in C_b^2(\mbb{R}_+)$ is a {\rv nonnegative, nondecreasing and concave} function and $\phi\in C_b^2(\mbb{R}_+)$ is a nonnegative nonincreasing function. Let
 \beqlb\label{eq4.6}
F(x,y)= \phi(x)f(x-y)\I_{\{x\neq y\}}, \quad (x,y)\in D.
 \eeqlb
Then $F\in \mcr{D}(\tilde{L})$ and, for $(x,y)\in \itDelta^c$,
 \beqlb\label{eq4.7}
\tilde{L}F(x,y)\ar\le\ar cy\phi(x) f''(x-y) + \big[\beta-bx-g(x)\big]\phi'(x)f(x-y) \\
 \ar\ar
+\,y\phi(x)\big[f(2(x-y))-2f(x-y)\big] \mu_{x-y}(0,\infty) \nnm\\
 \ar\ar
+\,\phi(x)\big[f(2(x-y))-2f(x-y)\big]\nu_{x-y}(0,\infty) \nnm\\
 \ar\ar
+\,\tilde{A}_1F(x,y) + \tilde{A}_2F(x,y), \nnm
 \eeqlb
where
 \beqlb\label{eq4.8}
\tilde{A}_1F(x,y)\ar=\ar f(x-y)\bigg[cx\phi''(x) - y\int_0^\infty \big[\phi(x+z)-\phi(x)\big] \mu_{y-x}(\d z) \\
 \ar\ar\qqquad
+\,x\int_0^\infty \big[\phi(x+z)-\phi(x)-\phi'(x)z\I_{\{z\le 1\}}\big] \mu(\d z) \nnm\\
 \ar\ar\qqquad
+\,\int_0^\infty \big[\phi(x+z)-\phi(x)\big] (\nu-\nu_{y-x})(\d z)\bigg] \nnm
 \eeqlb
and
 \beqlb\label{eq4.9}
\tilde{A}_2F(x,y)\ar=\ar (x-y)\phi(x)\bigg[cf''(x-y) - bf'(x-y) + \int_0^\infty \big[f(x-y+z) \\
 \ar\ar\qqquad
-\,f(x-y)-f'(x-y)z\I_{\{z\le 1\}}\big] \mu(\d z)\bigg]. \nnm
 \eeqlb
\elemma

\begin{proof} According to the definition of $F$, we have $F(z,z)= 0$ for $z\ge 0$. By \eqref{eq4.2} and \eqref{eq4.3} one can see that $F\in \mcr{D}(\tilde{L})$ and, for $(x,y)\in \itDelta^c$,
 \beqnn
\tilde{L}_0F(x,y)\ar=\ar \big[\beta-bx-g(x)\big]\phi'(x)f(x-y) - \big[bx-by+g(x)-g(y)\big]\phi(x)f'(x-y) \nnm\\
 \ar\ar
+\, c(x+3y)\phi(x)f''(x-y) + 2c(x+y)\phi'(x)f'(x-y) + cx\phi''(x)f(x-y) \nnm\\
 \ar\ar
+\,(x-y)\int_0^\infty \big[\phi(x+z)f(x-y+z)-\phi(x)f(x-y) \nnm\\
 \ar\ar\qqquad
-\,\phi(x)f'(x-y)z\I_{\{z\le 1\}} - \phi'(x)f(x-y)z\I_{\{z\le 1\}}\big] \mu(\d z) \nnm\\
 \ar\ar
+\,yf(x-y)\int_0^\infty \big[\phi(x+z)-\phi(x)-\phi'(x)z\I_{\{z\le 1\}}\big] \mu(\d z) \cr
 \ar\ar
+\,f(x-y)\int_0^\infty \big[\phi(x+z)-\phi(x)\big] \nu(\d z)
 \eeqnn
and
 \beqnn
\tilde{L}_1F(x,y)\ar=\ar y\big[f(2(x-y))-f(x-y)\big]\int_0^\infty \phi(x+z) \mu_{x-y}(\d z) \cr
 \ar\ar
-\, yf(x-y)\int_0^\infty \phi(x+z) \mu_{y-x}(\d z) - f(x-y)\int_0^\infty \phi(x+z) \nu_{y-x}(\d z) \cr
 \ar\ar
+\,\big[f(2(x-y))-f(x-y)\big]\int_0^\infty \phi(x+z) \nu_{x-y}(\d z).
 \eeqnn
Recall that $f$ and $g$ are nondecreasing and $\phi$ is nonincreasing. Then we have
 \beqnn
\tilde{L}_0F(x,y)\ar\le\ar \big[\beta-bx-g(x)\big]\phi'(x)f(x-y) - b(x-y)\phi(x)f'(x-y) \nnm\\
 \ar\ar
+\, cx\phi(x)f''(x-y) + cx\phi''(x)f(x-y) \nnm\\
 \ar\ar
+\,(x-y)\int_0^\infty \big[\phi(x+z)f(x-y+z) - \phi(x)f(x-y+z)\cr
 \ar\ar\qqquad\qqquad\qqquad
-\,\phi'(x)f(x-y)z\I_{\{z\le 1\}}\big] \mu(\d z) \cr
 \ar\ar
+\,(x-y)\phi(x)\int_0^\infty \big[f(x-y+z)-f(x-y)-f'(x-y)z\I_{\{z\le 1\}}\big] \mu(\d z) \nnm\\
\ar\ar
+\,yf(x-y)\int_0^\infty \big[\phi(x+z)-\phi(x)-\phi'(x)z\I_{\{z\le 1\}}\big] \mu(\d z) \cr
 \ar\ar
+\,f(x-y)\int_0^\infty \big[\phi(x+z)-\phi(x)\big] \nu(\d z) \nnm\\
 \ar\le\ar
\big[\beta-bx-g(x)\big]\phi'(x)f(x-y) - b(x-y)\phi(x)f'(x-y) \nnm\\
 \ar\ar
+\,cx\phi(x)f''(x-y) + cx\phi''(x)f(x-y) \nnm\\
 \ar\ar
+\,xf(x-y)\int_0^\infty \big[\phi(x+z)-\phi(x)-\phi'(x)z\I_{\{z\le 1\}}\big] \mu(\d z) \cr
 \ar\ar
+\,(x-y)\phi(x)\int_0^\infty \big[f(x-y+z)-f(x-y)-f'(x-y)z\I_{\{z\le 1\}}\big] \mu(\d z) \cr
 \ar\ar
+\,f(x-y)\int_0^\infty \big[\phi(x+z)-\phi(x)\big] \nu(\d z)
 \eeqnn
and
 \beqnn
\tilde{L}_1F(x,y)\ar=\ar y\phi(x)\big[f(2(x-y))-2f(x-y)\big]\mu_{x-y}(0,\infty) \nnm\\
 \ar\ar
+\, y\big[f(2(x-y))-f(x-y)\big]\int_0^\infty \big[\phi(x+z)-\phi(x)\big] \mu_{x-y}(\d z) \cr
 \ar\ar
-\, yf(x-y)\int_0^\infty \big[\phi(x+z)-\phi(x)\big] \mu_{y-x}(\d z) \cr
 \ar\ar
+\,\phi(x)\big[f(2(x-y))-2f(x-y)\big]\nu_{x-y}(0,\infty) \nnm\\
 \ar\ar
+\,\big[f(2(x-y))-f(x-y)\big]\int_0^\infty \big[\phi(x+z)-\phi(x)\big] \nu_{x-y}(\d z) \cr
 \ar\ar
-\, f(x-y)\int_0^\infty \big[\phi(x+z)-\phi(x)\big] \nu_{y-x}(\d z) \nnm\\
 \ar\le\ar
y\phi(x)\big[f(2(x-y))-2f(x-y)\big]\mu_{x-y}(0,\infty) \nnm\\
 \ar\ar
-\, yf(x-y)\int_0^\infty \big[\phi(x+z)-\phi(x)\big] \mu_{y-x}(\d z) \cr
 \ar\ar
+\,\phi(x)\big[f(2(x-y))-2f(x-y)\big]\nu_{x-y}(0,\infty) \nnm\\
 \ar\ar
-f(x-y)\int_0^\infty \big[\phi(x+z)-\phi(x)\big] \nu_{y-x}(\d z),
 \eeqnn
where we have used the relation $\mu_{x-y}(0,\infty)= \mu_{y-x}(0,\infty)$. Returning to \eqref{eq4.1} and reorganizing the terms we get \eqref{eq4.7}. \end{proof}

{\rv

In the sequel, we assume there is some $\lambda_0>0$ such that $\itPsi(\lambda_0)>0$. Under Grey's condition \eqref{eq1.9}, we can define some $F_0\in \mcr{D}(\tilde{L})$ by choosing an explicit form of the function $f\in C_b^2(\mbb{R}_+)$ in \eqref{eq4.6'}. To do so, take any $\lambda_0> 0$ such that $\itPsi(\lambda)> 0$ for $\lambda\ge \lambda_0$ and
 \beqnn
\int_{\lambda_0}^\infty \frac{1}{\itPsi(\lambda)} \d\lambda\le 1.
 \eeqnn
Fix a constant $\theta\ge 4$ and define
 \beqlb\label{eq4.11b}
F_0(x,y)= \big[\theta+h(x-y)\big]\I_{\{x\neq y\}}, \quad (x,y)\in D,
 \eeqlb
where
 \beqnn
h(z)= \int_{\lambda_0}^\infty \frac{1-\e^{-\lambda z}}{\itPsi(\lambda)} \d\lambda.
 \eeqnn
By Lemma~\ref{th4.2'} we have $F_0\in \mcr{D}(\tilde{L})$. It is easy to see that
 \beqlb\label{eq4.13'}
\theta\le F_0(x,y)\le 1+\theta, \quad (x,y)\in \itDelta^c.
 \eeqlb

\bproposition\label{th4.4'} Suppose that Condition~\ref{cd1.2}-(i) is satisfied. Let $F_0\in \mcr{D}(\tilde{L})$ be defined by \eqref{eq4.11b}. Then for any $l\ge 1$ there is a constant $\lambda_2> 0$ such that
 \beqlb\label{eq4.17'}
\tilde{L}F_0(x,y)\le -\lambda_2 F_0(x,y)\I_{\{x-y\le l\}}, \quad (x,y)\in \itDelta^c.
 \eeqlb
\eproposition

\begin{proof} It is clear that $h\in C_b^2(\mbb{R}_+)$ and the function is a nonnegative, nondecreasing and concave. By \eqref{eq4.7'} one can see that, for $(x,y)\in \itDelta^c$,
 \beqnn
\tilde{L}F_0(x,y)\ar\le\ar -\int_{\lambda_0}^\infty (x-y)\e^{-\lambda(x-y)} \d\lambda
 =
-\e^{-\lambda_0(x-y)}.
 \eeqnn
For $l\ge 1$ let $\lambda_2= (1+\theta)^{-1}{\e^{-\theta l}}$. Then $\tilde{L}F_0(x,y)\le -\e^{-\theta l}\le -\lambda_2 F_0(x,y)$ when $0< x-y\le l$. \end{proof}

}

{\rv We may also define a function $F_0\in \mcr{D}(\tilde{L})$ by special choices of the functions $f\in C_b^2(\mbb{R}_+)$ and $\phi\in C_b^2(\mbb{R}_+)$ in \eqref{eq4.6}. To do so, fix $\lambda_0>0$} such that $\itPsi(\lambda_0)>0$ and let
 \beqlb\label{eq4.10}
\psi(x)= 1-\e^{-\lambda_0x}, \quad x\ge 0.
 \eeqlb
For constants $0< x_0< 1$ and $\theta\ge 4$ to be specified later, we define
 \beqlb\label{eq4.11}
\phi(x)=\bigg\{\begin{array}{ll}
 \theta+(1-x/x_0)^3, &\quad 0\le x< x_0, \smallskip\cr
 \theta, &\quad x\ge x_0.
\end{array}
 \eeqlb
Then define the function
 \beqlb\label{eq4.12}
F_0(x,y)= \phi(x)\big[1+\psi(x-y)\big]\I_{\{x\neq y\}}, \quad (x,y)\in D.
 \eeqlb
By Lemma~\ref{th4.2} we have $F_0\in \mcr{D}(\tilde{L})$. It is easy to see that
 \beqlb\label{eq4.13}
\theta\le F_0(x,y)\le 2(1+\theta), \quad (x,y)\in \itDelta^c.
 \eeqlb

\blemma\label{th4.3} {\rv Let $F_0\in \mcr{D}(\tilde{L})$ be a function of the form \eqref{eq4.12}. Then for} any $l\ge 1$ there is a constant $\lambda_1> 0$ such that, for $(x,y)\in \itDelta^c$,
 \beqlb\label{eq4.14}
~~\qquad \tilde{L}F_0(x,y)\ar\le\ar - c\lambda_0^2\theta y\e^{-\lambda_0(x-y)} - \theta y\mu_{x-y}(0,\infty) - \lambda_1\theta\psi(x-y)\I_{\{x-y\le l\}} \\
 \ar\ar
-\, \theta\nu_{x-y}(0,\infty) + \big[y\mu_{x-y}(0,\infty) + I(x) + J(x)\big]\big[1+\psi(x-y)\big]\I_{\{x\le x_0\}} \nnm\\
 \ar\ar
+\,\big[1+\psi(x-y)\big]\int_0^\infty \big[\phi(x+z)-\phi(x)\big] (\nu-\nu_{y-x})(\d z)\I_{\{x\le x_0\}}, \nnm
 \eeqlb
where
 \beqlb\label{eq4.15}
I(x)= [\beta-bx-g(x)]\phi'(x)
 \eeqlb
and
 \beqlb\label{eq4.16}
J(x)= \frac{3x}{x_0^2}\bigg(2c + \int_0^1 z^2 \mu(\d z)\bigg).
 \eeqlb
\elemma

\begin{proof} Since $\psi(0)= 0$ and $\psi$ is concave, we have $\psi(2z)\le 2\psi(z)$ for $z\ge 0$. Then from \eqref{eq4.7} it follows that, for $(x,y)\in \itDelta^c$,
 \beqnn
\tilde{L}F_0(x,y)\ar\le\ar cy\phi(x) \psi''(x-y) + [\beta-bx-g(x)]\phi'(x)\big[1+\psi(x-y)\big] \nnm\\
 \ar\ar
+\,y\phi(x)\big[\psi(2(x-y))-2\psi(x-y)-1\big] \mu_{x-y}(0,\infty) \nnm\\
 \ar\ar
+\,\phi(x)\big[\psi(2(x-y))-2\psi(x-y)-1\big] \nu_{x-y}(0,\infty) \nnm\\
 \ar\ar
+\, \tilde{A}_1F_0(x,y) + \tilde{A}_2F_0(x,y) \nnm\\
 \ar\le\ar
-c\lambda_0^2\theta y\e^{-\lambda_0(x-y)} + I(x)\big[1+\psi(x-y)\big] - \theta y\mu_{x-y}(0,\infty) \nnm\\
 \ar\ar
-\,\theta\nu_{x-y}(0,\infty) + \tilde{A}_1F_0(x,y) + \tilde{A}_2F_0(x,y),
 \eeqnn
where we have also used the fact that $\phi(x)\ge \theta$. By the definition of $\phi$, we know that for all $x,z\ge 0$,
 \beqnn
\ar\ar 0\le \phi(x)-\phi(x+z)\le \I_{\{x\le x_0\}},
 \quad
0\le -\phi'(x)z\le \frac{3z}{x_0}\I_{\{x\le x_0\}}, \cr
\ar\ar \phi''(x)\le \frac{6}{x_0^2}\I_{\{x\le x_0\}},
 \quad
\phi(x+z)-\phi(x)-\phi'(x)z\le \frac{3z^2}{x_0^2}\I_{\{x\le x_0\}}.
 \eeqnn
From \eqref{eq4.8} it follows that
 \beqnn
\tilde{A}_1F_0(x,y)\ar=\ar \big[1+\psi(x-y)\big]\bigg[cx\phi''(x) - y\int_0^\infty \big[\phi(x+z)-\phi(x)\big] \mu_{y-x}(\d z) \cr
 \ar\ar\qqquad\qquad
+\,x\int_0^1 \big[\phi(x+z)-\phi(x)-\phi'(x)z\big] \mu(\d z) \cr
 \ar\ar\qqquad\qquad
+ \int_0^\infty \big[\phi(x+z)-\phi(x)\big] (\nu-\nu_{x-y})(\d z)\bigg] \cr
 \ar\le\ar
\big[1+\psi(x-y)\big]\bigg[\frac{6cx}{x_0^2} + y\mu_{x-y}(0,\infty) + x\int_0^1 \frac{3z^2}{x_0^2} \mu(\d z)\bigg]\I_{\{x\le x_0\}} \cr
 \ar\ar
+\,\big[1+\psi(x-y)\big]\int_0^\infty \big[\phi(x+z)-\phi(x)\big] (\nu-\nu_{y-x})(\d z)\I_{\{x\le x_0\}} \cr
 \ar=\ar
\big[1+\psi(x-y)\big]\big[y\mu_{x-y}(0,\infty) + J(x)\big]\I_{\{x\le x_0\}} \cr
 \ar\ar
+\,\big[1+\psi(x-y)\big]\int_0^\infty \big[\phi(x+z)-\phi(x)\big] (\nu-\nu_{y-x})(\d z)\I_{\{x\le x_0\}}.
 \eeqnn
In view of \eqref{eq4.9}, we have
 \beqnn
\tilde{A}_2F_0(x,y)\ar=\ar (x-y)\phi(x)\bigg[c\psi''(x-y) - b\psi'(x-y) + \int_0^\infty \big[\psi(x+z-y) \nnm\\
 \ar\ar\qqquad\qquad
-\, \psi(x-y) - \psi'(x-y)z\I_{\{z\le 1\}}\big] \mu(\d z)\bigg] \nnm\\
 \ar=\ar
- (x-y)\phi(x) \e^{-\lambda_0(x-y)}\itPsi(\lambda_0) \nnm\\
 \ar\le\ar
-\, \lambda_0^{-1}(1-\e^{-\lambda_0(x-y)})\theta \e^{-\lambda_0l}\itPsi(\lambda_0)\I_{\{x-y\le l\}}.
 \eeqnn
By putting together the above estimates we get \eqref{eq4.14} with $\lambda_1= \lambda_0^{-1}\e^{-\lambda_0l}\itPsi(\lambda_0)> 0$. \end{proof}

\bproposition\label{th4.4} {\rv Suppose that Condition~\ref{cd1.2}-(ii) is satisfied. Let $F_0\in \mcr{D}(\tilde{L})$ be a function defined by \eqref{eq4.12}. Then for any} $l\ge 1$ there is a constant $\lambda_2> 0$ such that
 \beqlb\label{eq4.17}
\tilde{L}F_0(x,y)\le -\lambda_2 F_0(x,y)\I_{\{x-y\le l\}}, \quad (x,y)\in \itDelta^c.
 \eeqlb
\eproposition

\begin{proof} The idea of the proof is to identify and take the advantage of the dominating factor among the branching, immigration and competition mechanisms in different parts of the space $\itDelta^c$. Under Condition~\ref{cd1.2}-(ii), we can choose constants $\kappa>0$ and $x_0\in (0,c_0)$ such that
 \beqlb\label{eq4.18}
c\lambda_0^2\e^{-\lambda_0x}+\mu_x(0,\infty)+\nu_x(0,\infty)\ge 2\kappa, \quad 0\le x\le x_0.
 \eeqlb
Note that $I(x)=0$ for all $x>x_0$. By \eqref{eq4.14} we have, for $x> x_0$,
 \beqlb\label{eq4.19}
\tilde{L}F_0(x,y)\ar\le\ar -\theta y\big[c\lambda_0^2\e^{-\lambda_0(x-y)}+\mu_{x-y}(0,\infty)\big] - \theta\nu_{x-y}(0,\infty) \\
 \ar\ar
-\, \lambda_1\theta\psi(x-y)\I_{\{x-y\le l\}}. \nnm
 \eeqlb
Since $x_0\le 1\le l$, using \eqref{eq4.14} again we see, for $0\le x\le x_0$,
 \beqlb\label{eq4.20}
\tilde{L}F_0(x,y)\ar\le\ar -c\lambda_0^2\theta y\e^{-\lambda_0(x-y)} -y\mu_{x-y}(0,\infty)\big[\theta-1-\psi(x-y)\big] \\
 \ar\ar
-\,\theta\nu_{x-y}(0,\infty) - \lambda_1\theta\psi(x-y) + \big[I(x)+J(x)\big]\big[1+\psi(x-y)\big] \nnm\\
 \ar\ar
+\,\big[1+\psi(x-y)\big]\int_0^\infty \big[\phi(x+z)-\phi(x)\big] (\nu-\nu_{y-x})(\d z). \nnm
 \eeqlb
With these estimates at hand, we prove the desired assertion by considering the following three cases.

{\rm(i)} We first consider the case of $x> x_0$. We will apply \eqref{eq4.19} and consider separately three subcases. From \eqref{eq4.19} it is easy to see $\tilde{L}F_0(x,y)\le 0$ when $x-y>l$. When $x_0/2< x-y\le l$, noticing that $\psi$ is nondecreasing and using the facts $\phi\le 1+\theta$ and $\psi\le 1$, we have
 \beqnn
\tilde{L}F_0(x,y)\ar\le\ar -\lambda_1\theta\psi(x-y)\le -\frac{\lambda_1\theta}{2}[\psi(x-y)+\psi(x_0/2)] \cr
 \ar\le\ar
-\frac{\lambda_1\theta}{2}\psi(x_0/2) \big[1+\psi(x-y)\big]
 \le
-\frac{\lambda_1\theta\psi(x_0/2)}{2(1+\theta)}F_0(x,y).
 \eeqnn
When $x-y\le x_0/2$, we have $y\ge x_0/2$ and so, by \eqref{eq4.18} and \eqref{eq4.19},
 \beqnn
\tilde{L}F_0(x,y)\ar\le\ar -\frac{\theta x_0}{2}\big[c\lambda_0^2\e^{-\lambda_0(x-y)} + \mu_{x-y}(0,\infty) + \nu_{x-y}(0,\infty)\big] - \lambda_1\theta\psi(x-y) \cr
 \ar\le\ar
-\theta x_0\kappa-\lambda_1\theta\psi(x-y)
 \le
-(x_0\kappa\land\lambda_1)\theta\big[1+\psi(x-y)\big] \cr
 \ar\le\ar
-(x_0\kappa\land\lambda_1)\frac{\theta}{1+\theta} F_0(x,y).
 \eeqnn

{\rm(ii)} Let us consider the case of $x\in(0,rx_0]$, where $r\in (0,1/2]$ will be specified later. Recall that $\theta\ge 4$ and $\psi\le 1$. Then, according to \eqref{eq4.20}, for $x\in [0,rx_0]$,
 \beqnn
\tilde{L}F_0(x,y)\ar\le\ar \big[I(x)+J(x)\big]\big[1+\psi(x-y)\big] - \theta\nu_{x-y}(0,\infty) \nnm\\
 \ar\ar
+\,\big[1+\psi(x-y)\big]\int_0^\infty \big[\phi(x+z)-\phi(x)\big] (\nu-\nu_{y-x})(\d z)\I_{\{x\le x_0\}} \nnm\\
 \ar\le\ar
\big[1+\psi(x-y)\big]\bigg[I(x)+J(x)-\nu_{x-y}(0,\infty) \nnm\\
 \ar\ar
+\int_0^\infty \big[\phi(x+z)-\phi(x)\big](\nu-\nu_{y-x})(\d z)\bigg].
 \eeqnn
Since $0< x_0< 1$, we can choose $r_*\in (0,1/2]$ and $q>0$ such that, for $0\le x\le r_*x_0$,
 \beqlb\label{eq4.21}
~\quad\ar\ar I(x) - \nu_{x-y}(0,\infty) + \int_0^\infty \big[\phi(x+z)-\phi(x)\big](\nu-\nu_{y-x})(\d z) \\
 \ar\ar\qquad
\le I(x) - \nu_{x-y}(0,\infty) - \int_{x_0-x}^\infty \Big(1-\frac{x}{x_0}\Big)^3 (\nu-\nu_{y-x})(\d z) \nnm\\
 \ar\ar\qqquad
- \int^{x_0-x}_0 \Big[\Big(1-\frac{x}{x_0}\Big)^3 - \Big(1-\frac{x+z}{x_0}\Big)^3\Big] (\nu-\nu_{y-x})(\d z) \nnm\\
 \ar\ar\qquad
\le\frac{3}{x_0}\big[|b|x+g(x)-\beta\big]\Big(1-\frac{x}{x_0}\Big)^2 - \nu_{x-y}(0,\infty) - \frac{1}{8}\int_{x_0-x}^\infty (\nu-\nu_{y-x})(\d z) \nnm\\
 \ar\ar\qqquad
- \frac{1}{x_0^3}\int^{x_0-x}_0 z^3 (\nu-\nu_{y-x})(\d z) \nnm\\
 \ar\ar\qquad
\le\frac{3}{x_0}\big[|b|x+g(x)\big] - \frac{3\beta}{4x_0} - \frac{1}{8}\int_0^\infty (1\land z^3) \nu(\d z)\le -q, \nnm
 \eeqlb
where for the last inequality we have used the {\rv Condition~\ref{cd1.1} for $\itPhi$ and} the fact that $g$ is a continuous function with $g(0)=0$. Now we take
 \beqlb\label{eq4.22}
r= r_*\land\Big[\frac{x_0q}{6}\Big(2c+\int_0^1 z^2 \mu(\d z)\Big)^{-1}\Big].
 \eeqlb
From \eqref{eq4.16} it follows that, for $x\in [0,rx_0]$,
 \beqnn
J(x)\le \frac{3r}{x_0}\Big(2c+\int_0^1 z^2 \mu(\d z)\Big)\le \frac{q}{2}.
 \eeqnn
Recall that $\theta\ge 4$ and $\psi\le 1$. Then, for $x\in [0,rx_0]$,
 \beqnn
\tilde{L}F_0(x,y)\ar\le\ar -\big[q-J(x)\big]\big[1+\psi(x-y)\big] \nnm\\
 \ar\le\ar
-\frac{q}{2}\big[1+\psi(x-y)\big]\le -\frac{q}{2(1+\theta)}F_0(x,y).
 \eeqnn

{\rm(iii)} We finally consider the case of $x\in(rx_0, x_0]$. In this case, we have
 \beqnn
I(x)\le \frac{3}{x_0}\big[|b|x_0+g(x_0)\big],
 \quad
J(x)\le \frac{3}{x_0}\bigg(2c+\int_0^1 z^2 \mu(\d z)\bigg).
 \eeqnn
Now we take
 \beqlb\label{eq4.23}
H= \frac{3}{x_0}\bigg[2c+|b|x_0+g(x_0)+\int_0^1 z^2 \mu(\d z)\bigg]
 \eeqlb
and
 \beqlb\label{eq4.24}
\theta=\max\Big\{4,\frac{2H}{\lambda_1},\frac{4H}{r\kappa x_0}, \frac{8H}{\lambda_1 \psi(rx_0/2)}\Big\}.
 \eeqlb
From \eqref{eq4.20} it follows that
 \beqlb\label{eq4.25}
~\qquad \tilde{L}F_0(x,y)\ar\le\ar -c\lambda_0^2\theta y\e^{-\lambda_0(x-y)} -y\mu_{x-y}(0,\infty) \big[\theta-1-\psi(x-y)\big] -\theta\nu_{x-y}(0,\infty) \\
 \ar\ar
+\,\psi(x-y)\big[I(x)+J(x)-\lambda_1\theta\big] + I(x) + J(x) \nnm\\
 \ar\le\ar
-\,c\lambda_0^2\theta y\e^{-\lambda_0(x-y)}-y\mu_{x-y}(0,\infty)\big[\theta-1-\psi(x_0)\big] -\theta\nu_{x-y}(0,\infty) \nnm\\
 \ar\ar
+\,\psi(x-y)(H-\lambda_1\theta) + H \nnm\\
 \ar\le\ar
-\,c\lambda_0^2\theta y\e^{-\lambda_0(x-y)} -\frac{\theta}{2}y\mu_{x-y}(0,\infty) -\theta\nu_{x-y}(0,\infty) \nnm\\
 \ar\ar
-\,\frac{\lambda_1\theta}{2} \psi(x-y) + H. \nnm
 \eeqlb
When $y\ge rx_0/2$, we have $x-y\le rx_0/2$ and, by \eqref{eq4.18},
 \beqnn
-c\lambda_0^2\theta y\e^{-\lambda_0(x-y)} -\frac{\theta}{2}y\mu_{x-y}(0,\infty) -\theta\nu_{x-y}(0,\infty)
 \le
-\frac{\theta r\kappa x_0}{2}.
 \eeqnn
Returning to \eqref{eq4.25}, we obtain
 \beqnn
\tilde{L}F_0(x,y)\ar\le\ar-\frac{\theta r\kappa x_0}{2}-\frac{\lambda_1\theta}{2}\psi(x-y)+H \cr
 \ar\le\ar
- \frac{\theta r\kappa x_0}{4}-\frac{\lambda_1\theta}{2}\psi(x-y) \cr
 \ar\le\ar
- \Big(\frac{r\kappa x_0}{2}\land\lambda_1\Big)\frac{\theta}{2}\big[1+\psi(x-y)\big] \cr
 \ar\le\ar
- \Big(\frac{r\kappa x_0}{2}\land\lambda_1\Big)\frac{\theta}{2(1+\theta)} F_0(x,y).
 \eeqnn
When $y\in [0,rx_0/2)$, we have $x-y\ge rx_0/2$ and hence
 \beqnn
\tilde{L}F_0(x,y)\ar\le\ar - \frac{\lambda_1\theta}{2}\psi(x-y) + H \cr
 \ar\le\ar
- \frac{\lambda_1\theta}{4} \psi(x-y) - \frac{\lambda_1\theta}{4}\psi(rx_0/2) + H \cr
 \ar\le\ar
- \frac{\lambda_1\theta}{8}\psi(x-y) - \frac{\lambda_1\theta}{8}\psi(rx_0/2) \cr
 \ar\le\ar
- \frac{\lambda_1\theta\psi(rx_0/2)}{8}\big[1+\psi(x-y)\big] \cr
 \ar\le\ar
- \frac{\lambda_1\theta\psi(rx_0/2)}{8(1+\theta)}F_0(x,y).
 \eeqnn
Combining all the estimates in the three cases, we obtain \eqref{eq4.17}. \end{proof}

\section{The exponential ergodicity}\label{section3}

 \setcounter{equation}{0}

In this section, we prove the exponential contraction property \eqref{eq1.31} for a suitable control function. From this property we derive the exponential ergodicity of the CBIC-process. {\rv Throughout the section, we assume that Conditions~\ref{cd1.1}, \ref{cd1.2} and~\ref{cd1.3} are satisfied. By Proposition~\ref{th2.3}, the CBIC-process with an arbitrary initial value is conservative. Let
 \beqlb\label{eq5.1}
V_2(x,y)= V(x)+ V(y),
 \quad
V_0(x,y)=V_2(x,y)\I_{\{x\neq y\}}, \quad (x,y)\in D.
 \eeqlb
Let $F_0\in \mcr{D}(\tilde{L})$ be given as in Propositions~\ref{th4.4'} and~\ref{th4.4} under Condition~\ref{cd1.2}-(i) and Condition~\ref{cd1.2}-(ii), respectively. For a constant $\varepsilon> 0$ to be specified later, define
 \beqlb\label{eq5.2}
G_0(x,y)= \varepsilon F_0(x,y) + V_0(x,y), \quad (x,y)\in D.
 \eeqlb
By \eqref{eq4.13'} and \eqref{eq4.13} it is easy to show that \eqref{eq1.32} holds.}

\bproposition\label{th5.1} {\rv We can define a function $G_0\in \mcr{D}(\tilde{L})$ by \eqref{eq5.2} such that, for some constant $\lambda_*>0$,}
 \beqlb\label{eq5.3}
\tilde{L}G_0(x,y)\le -\lambda_*G_0(x,y), \quad (x,y)\in \itDelta^c.
 \eeqlb
\eproposition

\begin{proof} Clearly, for any $(x,y)\in \itDelta^c$ the expression \eqref{eq4.2} of $\tilde{L}_0F(x,y)$ only depends on the restriction of $F$ to $\itDelta^c$. {\rv Then $V_0\in \mcr{D}(\tilde{L})$ and}
 \beqlb\label{eq5.4}
\tilde{L}_0V_0(x,y)= \tilde{L}_0V_2(x,y) = LV(x) + LV(y).
 \eeqlb
Since $0= V_0(z,z)\le V_2(z,z)$ for $z\ge 0$, by \eqref{eq4.3} we have $\tilde{L}_1V_0(x,y)\le \tilde{L}_1V_2(x,y)$. But, by \eqref{eq4.3}, \eqref{eq5.1} and the first equality in \eqref{eq3.2},
 \beqnn
\tilde{L}_1V_2(x,y)\ar=\ar y\int_0^\infty \big[V(2y+z-x)-V(y+z)\big] \mu_{x-y}(\d z) \cr
 \ar\ar
+\,y\int_0^\infty \big[V(x+z)-V(y+z)\big] \mu_{y-x}(\d z) \cr
 \ar\ar
+ \int_0^\infty \big[V(2y+z-x)-V(y+z)\big] \nu_{x-y}(\d z) \cr
 \ar\ar
+ \int_0^\infty \big[V(x+z)-V(y+z)\big] \nu_{y-x}(\d z)= 0.
 \eeqnn
It follows that
 \beqnn
\tilde{L}G_0(x,y)\le \varepsilon\tilde{L}F_0(x,y) + LV(x) + LV(y).
 \eeqnn
{\rv Let $l\ge 1$ be sufficiently large such that $V(z)> 12C_0/C_1$ for $z> l$ and let $\lambda_2> 0$ be the corresponding constant given by Propositions~\ref{th4.4'} and~\ref{th4.4} under Condition~\ref{cd1.2}-(i) and Condition~\ref{cd1.2}-(ii), respectively. Now let $\varepsilon= 4C_0/(\lambda_2\theta)$. By \eqref{eq1.20} we have}
 \beqnn
\tilde{L}G_0(x,y)\le -\varepsilon\lambda_2F_0(x,y)\I_{\{x-y\le l\}} + 2C_0 - C_1\big[V(x)+V(y)\big].
 \eeqnn
When $x\ge l$, {\rv since $\theta\ge 4$, we can use \eqref{eq4.13'} or \eqref{eq4.13}} to see
 \beqnn
\tilde{L} G_0(x,y)\ar\le\ar 2C_0 - C_1\big[V(x)+V(y)\big] \nnm\\
 \ar\le\ar
- 4C_0 - \frac{C_1}{2}\big[V(x)+V(y)\big] \nnm\\
 \ar\le\ar
-\frac{\lambda_2\theta}{2(1+\theta)}\varepsilon F_0(x,y) - \frac{C_1}{2}\big[V(x)+V(y)\big] \nnm\\
 \ar\le\ar\rv
-\frac{1}{2}\Big(\frac{4\lambda_2}{5}\land C_1\Big) G_0(x,y).
 \eeqnn
When $x\le l$, {\rv using \eqref{eq4.13'} or \eqref{eq4.13} again} we have
 \beqnn
\tilde{L} G_0(x,y)\ar\le\ar -\varepsilon\lambda_2 F_0(x,y) + 2C_0 - C_1\big[V(x)+V(y)\big] \nnm\\
 \ar\le\ar
-\frac{\varepsilon\lambda_2}{2} F_0(x,y) - \frac{\varepsilon\lambda_2\theta}{2} + 2C_0 - C_1\big[V(x)+V(y)\big] \nnm\\
 \ar\le\ar
-\frac{\varepsilon\lambda_2}{2} F_0(x,y) - C_1\big[V(x)+V(y)\big] \nnm\\
 \ar\le\ar
-\Big(\frac{\lambda_2}{2}\land C_1\Big) G_0(x,y).
 \eeqnn
Then \eqref{eq5.3} holds with $\rv \lambda_*= 2(C_1\land\lambda_2)/5> 0$. \end{proof}

\btheorem\label{th5.2} {\rv Let $G_0\in \mcr{D}(\tilde{L})$ and $\lambda_*> 0$ be given as in Proposition~\ref{th5.1}. Then we have \eqref{eq1.31} for $t\ge 0$ and $(x,y)\in D$.} \etheorem

\begin{proof} It suffices to consider $(x,y)\in \itDelta^c$. Let $\{(X_t,Y_t): t\ge 0\}$ be the Markov coupling defined by \eqref{eq3.8} and \eqref{eq3.9} with $(X_0,Y_0)= (x,y)$. Recall that $\zeta_n^*= \inf\{t\ge 0: X_t\ge n$ or $X_t-Y_t\le 1/n\}$. By Proposition~\ref{th2.3} and Theorem~\ref{th3.1}, the process $\{(X_t,Y_t): t\ge 0\}$ is conservative. Then \eqref{eq3.11} implies $\lim_{n\to \infty}\zeta_n^*= \lim_{n\to \infty}T_n= T$, where $T_n= \inf\{t\ge 0: X_t-Y_t\le 1/n\}$. By Theorem~\ref{th4.1} and integration by parts, for any $n\ge 1$ we have
 \beqlb\label{eq5.6}
~\quad \e^{\lambda_*(t\land \zeta_n^*)}G_0(X_{t\land \zeta_n^*},Y_{t\land \zeta_n^*})= G_0(x,y) + \int_0^{t\land \zeta_n^*} \e^{\lambda_*s}(\tilde{L}+\lambda)G_0(X_s,Y_s) \d s + M_n(t),
 \eeqlb
where $\{M_n(t): t\ge 0\}$ is a martingale. From \eqref{eq5.3} and \eqref{eq5.6} it follows that
 \beqnn
\mbf{E}\big[\e^{\lambda_*(t\land \zeta_n^*)}G_0(X_{t\land \zeta_n^*},Y_{t\land \zeta_n^*})\big]\le G_0(x,y),
 \eeqnn
and so
 \beqnn
\mbf{E}\big[\e^{\lambda_*(t\land \zeta_n^*)}G_0(X_{t\land \zeta_n^*},Y_{t\land \zeta_n^*})\I_{\{t< T\}}\big]\le G_0(x,y).
 \eeqnn
Since $G_0(x,x)= 0$ for $x\ge 0$ and $X_{T+t}= Y_{T+t}$ for $t\ge 0$, we can let $n\to \infty$ and use Fatou's lemma to get
 \beqnn
\mbf{E}\big[\e^{\lambda_*t}G_0(X_t,Y_t)\big]= \mbf{E}\big[\e^{\lambda_*t}G_0(X_t,Y_t)\I_{\{t< T\}}\big]\le G_0(x,y),
 \eeqnn
which clearly implies \eqref{eq1.31}. \end{proof}

The result of the next lemma should be already known, but we could not find a reference. For the convenience of the reader, we give a simple proof of the result here.

\blemma\label{th5.3} The expressions \eqref{eq1.14} and \eqref{eq1.16} for the $V$-weighted total variation distance $W_V$ are equivalent.
\elemma

\begin{proof} For $\gamma,\eta\in \mcr{P}_V(\mbb{R}_+)$ let $U= \supp((\gamma-\eta)_-)$ {\rv and $U^c= \mbb{R}_+\setminus U$, where $(\gamma-\eta)_-$ denotes the lower variation of the signed measure $\gamma-\eta$ in its Jordan decomposition.} If $\pi\in \mcr{C}(\gamma,\eta)$, then
 \beqnn
\int_{\mbb{R}_+^2} d_V(x,y) \pi(\d x,\d y)
 \ar=\ar
\int_{\mbb{R}_+^2} [2+V(x)+V(y)]\I_{\{x\neq y\}} \pi(\d x,\d y) \cr
 \ar\rv=\ar
\rv\int_{\mbb{R}_+^2} [2+V(x)+V(y)] \pi(\d x,\d y) - 2\int_{\mbb{R}_+^2} [1+V(x)]\I_{\{x=y\}} \pi(\d x,\d y) \cr
 \ar\rv=\ar
\int_{\mbb{R}_+} [1+V(x)] \gamma(\d x) + \int_{\mbb{R}_+} [1+V(y)] \eta(\d y) - 2\int_{U\times U} [1+V(x)] \pi(\d x,\d y) \cr
 \ar\ar
-\, 2\int_{U^c\times U^c} [1+V(y)] \pi(\d x,\d y) \cr
 \ar\ge\ar
\int_{\mbb{R}_+} [1+V(x)] \gamma(\d x) + \int_{\mbb{R}_+} [1+V(y)] \eta(\d y) - 2\int_{U} [1+V(x)] \gamma(\d x) \cr
 \ar\ar
-\, 2\int_{U^c} [1+V(y)] \eta(\d y) \cr
 \ar=\ar
\int_{\mbb{R}_+} [1+V(x)] \gamma(\d x) + \int_{\mbb{R}_+} [1+V(y)] \eta(\d y) - 2\int_{\mbb{R}_+} [1+V(x)] (\gamma\land \eta)(\d x) \cr
 \ar=\ar
\int_{\mbb{R}_+}[1+V(x)]|\gamma-\eta|(\d x).
 \eeqnn
One the other hand, let $\pi_*\in \mcr{C}(\gamma,\eta)$ be defined by
 \beqnn\rv
\pi_*(\d x,\d y)= (\gamma\land \eta)_*(\d x,\d y) + \frac{(\gamma-\eta)^+(\d x)(\gamma-\eta)^-(\d y)}{(\gamma-\eta)^+(\mbb{R}_+)}.
 \eeqnn
where $(\gamma\land \eta)_*$ is the image of $\gamma\land \eta$ under the mapping $x\mapsto (x,x)$ from $\mbb{R}_+$ to $\mbb{R}_+^2$. It is easy to see that
 \beqnn
\int_{\mbb{R}_+^2} d_V(x,y) \pi_*(\d x,\d y)
 \ar=\ar
\frac{1}{(\gamma-\eta)^+(\mbb{R}_+)}\int_{U\times U^c} [1+V(x)] (\gamma-\eta)^+(\d x)(\gamma-\eta)^-(\d y) \cr
 \ar\ar
+ \frac{1}{(\gamma-\eta)^+(\mbb{R}_+)}\int_{U\times U^c} [1+V(y)] (\gamma-\eta)^+(\d x)(\gamma-\eta)^-(\d y) \cr
 \ar=\ar
\int_U [1+V(x)] (\gamma-\eta)^+(\d x) + \int_{U^c} [1+V(y)] (\gamma-\eta)^-(\d y) \cr
 \ar=\ar
\int_{\mbb{R}_+}[1+V(x)]|\gamma-\eta|(\d x).
 \eeqnn
Those clearly imply the desired result. \end{proof}

\begin{proof}[Proof of Theorem~\ref{th1.1}] By Theorem~\ref{th5.2}, we have \eqref{eq1.31} for $(x,y)\in D$. It is easy to see that $\tilde{P}_t((x,y),\cdot)$ is a coupling of $P_t(x,\cdot)$ and $P_t(y,\cdot)$. Then \eqref{eq1.19} follows from \eqref{eq1.31} and \eqref{eq1.32} with $K= c_2/c_1$. By the convexity of the Wasserstein distance we have, for any $\gamma,\eta\in \mcr{P}_V(\mbb{R}_+)$ and $\pi\in \mcr{C}(\gamma,\eta)$,
 \beqnn
W_V(\gamma P_t,\eta P_t)
 \ar\le\ar
W_V\bigg(\int_{\mbb{R}_+^2}P_t(x,\cdot)\pi(\d x,\d y), \int_{\mbb{R}_+^2}P_t(y,\cdot)\pi(\d x,\d y)\bigg) \\
 \ar\le\ar
\int_{\mbb{R}_+^2} W_V(P_t(x,\cdot), P_t(y,\cdot))\pi(\d x,\d y)
 \le\rv
K\e^{-\lambda_*t} \int_{\mbb{R}_+^2} d_V(x,y)\pi(\d x,\d y);
 \eeqnn
see, e.g., Villani \cite[Theorem~4.8]{Vil09}. {\rv It follows that}
 \beqlb\label{eq5.7}\rv
W_V(\gamma P_t,\eta P_t)
 \le
K\e^{-\lambda_*t} \inf_{\pi\in \mcr{C}(\gamma,\eta)}\int_{\mbb{R}_+^2} d_V(x,y)\pi(\d x,\d y)
 =
K\e^{-\lambda_*t} W_V(\gamma,\eta).
 \eeqlb
Then for sufficiently large $r>0$, the operator $P_r$ on $\mcr{P}_V(\mbb{R}_+)$ is contractive. By the Banach fixed point theorem and the completeness of $\mcr{P}_V(\mbb{R}_+)$, there is a unique $\gamma_r\in \mcr{P}_V(\mbb{R}_+)$ such that $\gamma_r P_r=\gamma_r$. Now we fix such an $r>0$ and define $\gamma= r^{-1}\int_0^r\gamma_r P_s\d s$. By the Chapman-Kolmogorov equation, for $0\le t< r$ we have
\beqnn
\gamma P_t\ar=\ar \rv\frac{1}{r}\int_t^{r+t}\gamma_r P_s\d s
 =
\frac{1}{r}\int_t^r\gamma_r P_s\d s + \frac{1}{r}\int_0^t\gamma_r P_{r+s}\d s \cr
 \ar=\ar\rv
\frac{1}{r}\int_t^r\gamma_r P_s\d s + \frac{1}{r}\int_0^t\gamma_r P_s\d s
 =
\frac{1}{r}\int_0^r\gamma_r P_s\d s= \gamma.
 \eeqnn
More generally, for any $t\ge 0$ there is a unique integer $k\ge 0$ such that $0\le t-kr< r$. Then $\gamma P_t= \gamma P_{kr}P_{t-kr}= \gamma P_{t-kr}= \gamma$. By applying (\ref{eq5.7}) again we obtain (\ref{eq1.18}) with $C(\eta)= KW_V(\gamma,\eta)$. \end{proof}

\begin{proof}[Proofs of Propositions~\ref{th1.2} and~\ref{th1.3}] Under the integrability condition \eqref{eq1.22}, we have $V_1\in \mcr{D}(L)$. By \eqref{eq1.13} it is easy to see that
 \beqnn
LV_1(x)= [\beta-bx-g(x)] + x\int_1^\infty z \mu(\d z) + \int_0^\infty z \nu(\d z).
 \eeqnn
Then \eqref{eq1.20} is equivalent to \eqref{eq1.23}. The integrability condition \eqref{eq1.24} implies $V_{\log}\in \mcr{D}(L)$. By \eqref{eq1.13} we have
 \beqnn
LV_{\log}(x)\ar=\ar -\frac{c}{(1+x)^2} + x\int_0^\infty \Big[\log\Big(1+\frac{z}{1+x}\Big)-\frac{z}{1+x}\I_{\{z\le 1\}}\Big] \mu(\d z) \cr
 \ar\ar
+\, \frac{\beta-bx-g(x)}{1+x} + \int_0^\infty \log\Big(1 + \frac{z}{1+x}\Big) \nu(\d z).
 \eeqnn
By {\rv Taylor's expansion,} we have
 \beqnn
\frac{z}{1+x}-\log\Big(1+\frac{z}{1+x}\Big)
 =
z^2 \int_0^1\frac{(1-u)\d u}{(1+x+uz)^2} \cr
 \eeqnn
By the dominated convergence theorem,
 \beqnn
\lim_{x\to \infty} x^2\int_0^1 \Big[\frac{z}{1+x}-\log\Big(1+\frac{z}{1+x}\Big)\Big] \mu(\d z)= \frac{1}{2}\int_0^1 z^2 \mu(\d z).
 \eeqnn
It follows that
 \beqlb\label{eq5.8}
\lim_{x\to \infty} \frac{x}{\log x}\int_0^1 \Big[\log\Big(1+\frac{z}{1+x}\Big)-\frac{z}{1+x}\Big] \mu(\d z)=0.
 \eeqlb
Then \eqref{eq1.20} is equivalent to \eqref{eq1.25}. \end{proof}

\begin{proof}[Proofs of Corollaries~\ref{th1.4} and~\ref{th1.5}] Corollary~\ref{th1.4} follows easily by Proposition~\ref{th1.2}. Then it remains to prove Corollary~\ref{th1.5}. By the proof of Theorem~\ref{th2.4}, the CBIC-process with stable branching mechanism has generator $L$ defined by \eqref{eq1.13} with $b= a + \sigma h_\alpha$ and $\mu(\d z)= \alpha\sigma m_\alpha(\d z)$, where $m_\alpha$ and $h_\alpha$ are defined by \eqref{eq1.26} and \eqref{eq2.6}, respectively. For $\alpha= 1$, by elementary calculus we have
 \beqnn
\int_1^\infty \log\Big(1+\frac{z}{1+x}\Big) \frac{1}{z^2}\d z
 =
\Big(\frac{1}{1+x}\Big)\log(2+x) + \log\Big(\frac{2+x}{1+x}\Big).
 \eeqnn
For $0< \alpha< 1$, by Zwillinger \cite[3.194.4, p.318]{Zwi18} we have
 \beqnn
\int_0^\infty \log\Big(1+\frac{z}{1+x}\Big)\frac{1}{z^{1+\alpha}}\d z
 =
\frac{\pi}{\alpha\sin(\alpha\pi)(1+x)^\alpha}.
 \eeqnn
It follows that
 \beqnn
\int_1^\infty \log\Big(1+\frac{z}{1+x}\Big)\frac{1}{z^{1+\alpha}}\d z
 \ar=\ar
\frac{\pi}{\alpha\sin(\alpha\pi)(1+x)^\alpha} - \int_0^1 \log\Big(1+\frac{z}{1+x}\Big)\frac{1}{z^{1+\alpha}}\d z \cr
 \ar=\ar
\int_0^1 \Big[\frac{z}{1+x}-\log\Big(1+\frac{z}{1+x}\Big)\Big]\frac{1}{z^{1+\alpha}}\d z \cr
 \ar\ar
+ \frac{\pi}{\alpha\sin(\alpha\pi)(1+x)^\alpha} - \frac{1}{(1-\alpha)(1+x)}.
 \eeqnn
Then Corollary~\ref{th1.5} follows by Proposition~\ref{th1.3} and \eqref{eq5.8}. \end{proof}

\bremark\label{about_lambda} {\rv Our approach provides a way of finding the exponential ergodicity rate $\lambda_*> 0$. Let $C_0> 0$ and $C_1> 0$ be as in \eqref{eq1.20}. Then $\lambda_*= 2(C_1\land\lambda_2)/5$ by the proof of Proposition~\ref{th5.1}. {\rrv To determine $\lambda_2> 0$, under Condition~\ref{cd1.2}-(i) we follow the proof of Proposition~\ref{th4.4'}, while under Condition~\ref{cd1.2}-(ii) we follow the steps given below:}}
 \benumerate

\itm[(1)] determine $\kappa= \kappa(\lambda_0,c_0)$ and $x_0= x_0(\lambda_0,c_0)$ by \eqref{eq4.18};

\itm[(2)] choose $l= l(C_0,C_1)$ as in the proof of Proposition~\ref{th5.1};

\itm[(3)] let $\lambda_1= \lambda_0^{-1}\e^{-\lambda_0l}\itPsi(\lambda_0)$ as in the proof of Lemma~\ref{th4.3};

\itm[(4)] determine $q= q(x_0)$, $r_*= r_*(x_0)$ and $r= r(x_0)$ by \eqref{eq4.21} and \eqref{eq4.22};

\itm[(5)] define $H= H(x_0)$ and $\theta= \theta(\lambda_1,x_0,\kappa,r)$ by \eqref{eq4.23} and \eqref{eq4.24};

\itm[(6)] choose $\lambda_2= \lambda_2(\lambda_1,x_0,\kappa,r,\theta,q)$ as indicated in the proof of Proposition~\ref{th4.4}.

 \eenumerate
\eremark

{\rv

\bexample\label{ex5.1} Let $\itPsi(\lambda)= \lambda^\alpha$, $\itPhi(\lambda)= \lambda$ and $g(x)= x^2$, where $1< \alpha\le 2$. In this case, it is easy to see that Conditions~\ref{cd1.1} and~\ref{cd1.2}-(i) are satisfied and
 \beqnn
LV_1(x)= 1 - x^2\le 2 - V_1(x), \quad x\ge 0.
 \eeqnn
Then \eqref{eq1.20} holds for $V_1$ with $C_0= 2$ and $C_1= 1$. Take $\theta= 4$ in \eqref{eq4.11b}. According the proof of Proposition~\ref{th5.1}, we can choose $l= 12C_0/C_1= 24$. Next take $\lambda_2= (1+\theta)^{-1}{\e^{-\theta l}}= \e^{-96}/5$ as in the proof of Proposition~\ref{th4.4'}. Finally, we have $\lambda_*= 2(C_1\land\lambda_2)/5= 2\e^{-96}/25$. \qed \eexample

}

{\rv The estimates in the procedure {\rrv of determining the constant} $\lambda_*> 0$ are certainly not optimal. It {\rrv remains an interesting} problem to improve the arguments to get the optimal exponential ergodicity rate.}

 \begin{funding}
This research is supported by the National Key R{\&}D Program of China (Nos.~2020YFA0712901 {\rrv and 2022YFA1000033) and} the National Natural Science Foundation of China (Nos.~11731012, 11831014, 12271029, {\rrv 12071076 and 12225104).}
 \end{funding}


\begin{thebibliography}{4}\small
	

\bibitem{Ali85} Aliev, S.A. (1985): A limit theorem for the Galton--Watson branching processes with immigration. \textit{Ukrainian Math. J.} \textbf{37}, 535--438.

\bibitem{AlS82} Aliev, S.A. and Shchurenkov, V.M. (1982): Transitional phenomena and the convergence of Galton--Watson processes to Ji\v{r}ina processes. \textit{Theory Probab. Appl.} \textbf{27}, 472--485.	

\bibitem{AtN72} Athreya, K.B. and Ney, P.E. (1972): \textit{Branching Processes}. Springer, Heidelberg.

\bibitem{BBS13} Berestycki, J., Berestycki, N. and Schweinsberg, J. (2013): The genealogy of branching Brownian motion with absorption. \textit{Ann. Probab.} \textbf{41}, 527--618.

\bibitem{BFF18} Berestycki, J., Fittipaldi, M.C. and Fontbona, J. (2018): Ray-Knight representation of flows of branching processes with competition by pruning of L\'evy trees. \textit{Probab. Theory Related Fields} \textbf{172}, 725--788.

\bibitem{BeL00} Bertoin, J. and Le~Gall, J.-F. (2000): The Bolthausen--Sznitman coalescent and the genealogy of continuous-state branching processes. \textit{Probab. Theory Related Fields} \textbf{117}, 249--266.

\bibitem{BeL03} Bertoin, J. and Le~Gall, J.-F. (2003): Stochastic flows associated to coalescent processes. \textit{Probab. Theory Related Fields} \textbf{126}, 261--288.

\bibitem{BeL05} Bertoin, J. and Le~Gall, J.-F. (2005): Stochastic flows associated to coalescent processes II: Stochastic differential equations. \textit{Ann. Inst. H. Poincar\'e Probab. Statist.} \textbf{41}, 307--333.

\bibitem{BeL06} Bertoin, J. and Le~Gall, J.-F. (2006): Stochastic flows associated to coalescent processes III: Limit theorems. \textit{Illinois J. Math.} \textbf{50}, 147--181.

\bibitem{BoS98} Bolthausen, E. and Sznitman, A.-S. (1998): On Ruelle's probability cascades and an abstract cavity method. \textit{Commun. Math. Phys.} \textbf{197}, 247--276.

\bibitem{Che86a} {\rv Chen, M.-F. (1986a): \textit{Jump Processes and Interacting Particle Systems} (In Chinese). Beijing Normal Univ. Press, Beijing.}

\bibitem{Che86b} {\rv Chen, M.-F. (1986b): Couplings of jump processes. \textit{Acta Math. Sinica, New Series} \textbf{2}, 123--136.}

\bibitem{Che91} Chen, M.F. (1991): On three classical problems for Markov chains with continuous time parameters. \textit{J. Appl. Probab.} \textbf{28}, 305--320.

\bibitem{Che04} Chen, M.-F. (2004): \textit{From Markov Chains to Non-Equilibrium Particle Systems}. 2nd Ed. World Scientific, Singapore.

\bibitem{Che05} Chen, M.F. (2005): \textit{Eigenvalues, Inequality and Ergodic Theory}. Springer, London.

\bibitem{DaL06} Dawson, D.A. and Li, Z. (2006): Skew convolution semigroups and affine markov processes. \textit{Ann. Probab.} \textbf{34}, 1103--1142.

\bibitem{DaL12} Dawson, D.A. and Li, Z. (2012): Stochastic equations, flows and measure-valued processes. \textit{Ann. Probab.} \textbf{40}, 813--857.

 \bibitem{DMT95} Down, D., Meyn, S.P. and Tweedie, R.T. (1995): Exponential and uniform ergodicity of Markov processes. \textit{Ann. Probab.} \textbf{23}, 1671--1691.

\bibitem{Ebe11} Eberle, A. (2011): Reflection coupling and Wasserstein contractivity without convexity. \textit{C. R. Math. Acad. Sci. Paris} \textbf{349}, 1101--1104.

\bibitem{Ebe16} Eberle, A. (2016): Reflection couplings and contraction rates for diffusions. \textit{Probab. Theory Related Fields} \textbf{166}, 851--886.

\bibitem{EGZ19} Eberle, A., Guillin, A. and Zimmer, R. (2019): Quantitative Harris-type theorems for diffusions and McKean-Vlasov processes. \textit{Trans. Amer. Math. Soc.} \textbf{371}, 7135--7173.

\bibitem{Fel51} Feller, W. (1951): Diffusion processes in genetics. In: \textit{Proceedings {\rm2}nd Berkeley Symp. Math. Statist. Probab.} 1950, 227--246. Univ. of California Press, Berkeley and Los Angeles.

\bibitem{Fou19}{\rv Foucart, C. (2019): Continuous-state branching processes with competition: duality and reflection at infinity. \textit{Electron. J. Probab.} \textbf{24}, article no. 33, 1--38.}

\bibitem{Fri23} Friesen, M. (2023): Long-time behavior for subcritical measure-valued branching processes with immigration. \textit{Potential Anal.} \textbf{59}, 705--730.

\bibitem{FJKR23} Friesen, M., Jin, P., Kremer, J. and R\"{u}diger, B. (2023): Exponential ergodicity for stochastic equations of nonnegative processes with jumps. \textit{ALEA Lat. Am. J. Probab. Math. Stat.} \textbf{20}, 593--627.

\bibitem{FJR20} Friesen, M., Jin, P. and R\"{u}diger, B. (2020): Stochastic equation and exponential ergodicity in Wasserstein distances for affine processes. \textit{Ann. Appl. Probab.} \textbf{30}, 2165--2195.

\bibitem{FuL10} Fu, Z. and Li, Z. (2010): Stochastic equations of non-negative processes with jumps. \textit{Stochastic Process. Appl.} \textbf{120}, 306--330.

\bibitem{Gre74} Grey, D.R. (1974): Asymptotic behaviour of continuous time, continuous state-space branching processes. \textit{J. Appl. Probab.} \textbf{11}, 669--677.

\bibitem{Gri74} Grimvall, A. (1974): On the convergence of sequences of branching processes. \textit{Ann. Probab.} \textbf{2}, 1027--1045.

\bibitem{HaM11} Hairer, M. and Mattingly, J.C. (2011): Yet another look at Harris' ergodic theorem for Markov chains. In: \textit{Seminar on Stochastic Analysis, Random Fields and Applications} VI, Progr. Probab., Vol.~63, Birkh\"auser/Springer Basel AG, Basel, pp. 109--117.

\bibitem{Har63} Harris, T.E. (1963): \textit{The Theory of Branching Processes}. Springer, Heidelberg.

\bibitem{IkW89} Ikeda, N. and Watanabe, S. (1989): \textit{Stochastic Differential Equations and Diffusion Processes}. 2nd Ed. North-Holland, Amsterdam; Kodansha, Tokyo.

\bibitem{KaW71} Kawazu, K. and Watanabe, S. (1971): Branching processes with immigration and related limit theorems. \textit{Theory Probab. Appl.} \textbf{16}, 36--54.


\bibitem{Lam05} Lambert, A. (2005): The branching process with logistic growth. \textit{Ann. Appl. Probab.} \textbf{15}, 1506--1535.

\bibitem{Lam67a} Lamperti, J. (1967): The limit of a sequence of branching processes. \textit{Z. Wahrsch. Verw. Gebiete} \textbf{7}, 271--288.

\bibitem{Li19} Li, P.-S. (2019): A continuous-state polynomial branching process. \textit{Stochastic Process. Appl.}, \textbf{129}, 2941--2967.

\bibitem{LiW20} Li, P.-S. and Wang, J. (2020): Exponential ergodicity for general continuous-state nonlinear branching processes. \textit{Electron. J. Probab.} \textbf{25}, article no. 125, 1--25.

\bibitem{LYZ19} Li, P.-S., Yang, X. and Zhou, X. (2019): A general continuous-state nonlinear branching process. \textit{Ann. Appl. Probab.} \textbf{29} , 2523--2555.

\bibitem{Li22} {\rv Li, Z. (2022): \textit{Measure-Valued Branching Markov Processes}. 2nd Ed. Springer, Heidelberg.}

\bibitem{Li20} Li, Z. (2020): Continuous-state branching processes with immigration. A Chapter in: \textit{From Probability to Finance}, pp.\,1--69, edited by Y.~Jiao. \textit{Mathematical Lectures from Peking University}. Springer, Singapore.

\bibitem{Li21} Li, Z. (2021): Ergodicities and exponential ergodicities of Dawson-Watanabe type processes. \textit{Theory Probab. Appl.} \textbf{66}, 276-298.

\bibitem{LiM15} Li, Z. and Ma, C. (2015): Asymptotic properties of estimators in a stable Cox--Ingersoll--Ross model. \textit{Stochastic Process. Appl.} \textbf{125}, 3196--3233.

\bibitem{LMW21} Liang, M., Majka, M. and Wang, J. (2021): Exponential ergodicity for SDEs and McKean-Vlasov processes with L\'evy noise. \textit{Ann. Inst. Henri Poincar\'e Probab. Stat.} \textbf{57}, 1665--1701.

\bibitem{LuW16} Luo, D. and Wang, J. (2016): Exponential convergence in $L^p$-Wasserstein distance for diffusion processes without uniformly dissipative drift. \textit{Math. Nachr.} \textbf{289}, 1909--1926.

\bibitem{LuW19} Luo, D. and Wang, J. (2019): Refined basic couplings and Wasserstein-type distances for SDEs with L\'evy noises. \textit{Stochastic Process. Appl.} \textbf{129}, 3129--3173.

\bibitem{Maj17} Majka, M.B. (2017): Coupling and exponential ergodicity for stochastic differential equations driven by L\'evy processes. \textit{Stochastic Process. Appl.} \textbf{127}, 4083--4125.

\bibitem{MeT92} Meyn, S. and Tweedie, R.L. (1992): Stability of Markovian processes I: Criteria for discrete-time chains. \textit{Adv. Appl. Probab.} \textbf{24}, 542--574.

\bibitem{MeT93a} Meyn, S. and Tweedie, R.L. (1993): Stability of Markovian processes II: Continuous-time processes and sampled chains. \textit{Adv. Appl. Probab.} \textbf{25}, 487--517.

\bibitem{MeT93b} Meyn, S. and Tweedie, R.L. (1993): Stability of Markovian processes III: {\rv Foster-Lyapunov criteria} for continuous-time processes. \textit{Adv. Appl. Probab.} \textbf{25}, 518--548.

\bibitem{PaP18} Palau, S. and Pardo, J.C. (2018): Branching processes in a L\'evy random environment. \textit{Acta Appl. Math.} \textbf{153}, 55--79.

\bibitem{Par16} Pardoux, E. (2016): \textit{Probabilistic Models of Population Evolution: Scaling Limits, Genealogies and Interactions}. Springer, Switzerland.

\bibitem{Pin72} Pinsky, M.A. (1972): Limit theorems for continuous state branching processes with immigration. \textit{Bull. Amer. Math. Soc.} \textbf{78}, 242--244.

\bibitem{ScW12} Schilling, R.L. and Wang, J. (2012): On the coupling property and the Liouville theorem for Ornstein-Uhlenbeck processes. \textit{J. Evol. Equat.} \textbf{12}, 119-140.

\bibitem{Sit05} Situ, R. (2005): \textit{Theory of Stochastic Differential Equations with Jumps and Applications}. Springer, Heidelberg.

\bibitem{Sta03a} Stannat, W. (2003): Spectral properties for a class of continuous state branching processes with immigration. \textit{J. Funct. Anal.} \textbf{201}, 185--227.

\bibitem{Sta03b} Stannat, W. (2003): On transition semigroups of $(A,\itPsi)$-superprocesses with immigration. \textit{Ann. Probab.} \textbf{31}, 1377--1412.

\bibitem{Vil09} Villani, C. (2009): \textit{Optimal Transport, Old and New}. Springer, Berlin.

\bibitem{Zwi18} Zwillinger, D. (2018): \textit{Table of Integrals, Series, and Products}. 8th Edition, Elsevier, Sinagpore.


\end{thebibliography}
\end{document}